\documentclass[opre,nonblindrev,dvi-to-ps]{informs3}

\DoubleSpacedXI 

\usepackage{endnotes}
\let\footnote=\endnote
\let\enotesize=\normalsize
\def\notesname{Endnotes}%
\def\makeenmark{\hbox to1.275em{\theenmark.\enskip\hss}}
\def\enoteformat{\rightskip0pt\leftskip0pt\parindent=1.275em
  \leavevmode\llap{\makeenmark}}

\usepackage{xcolor}
\usepackage{amsfonts, amssymb}
\usepackage{algorithm, algorithmic}
\usepackage{tikz}
\usepackage{pgfplots}
\usepackage{subcaption}
\usepackage{ulem}
\usepackage[toc,page]{appendix}
\usepackage[justification=centering]{caption}

\newcommand{\bmu}{\pmb{\mu}}
\newcommand{\balpha}{\pmb{\alpha}}
\newcommand{\btheta}{\pmb{\theta}}
\newcommand{\blambda}{\pmb{\lambda}}
\newcommand{\bbeta}{\pmb{\beta}}
\newcommand{\bigeta}{\pmb{\eta}}

\newcommand{\bv}{{\bf v}}

\usepackage{mathtools}

\title{\fontsize{17}{20}\selectfont Consistency and Computation for Regularized Maximum Likelihood Estimation of Multivariate Hawkes Processes}

\usepackage{natbib}
 \bibpunct[, ]{(}{)}{,}{a}{}{,}%
 \def\bibfont{\small}%
\TheoremsNumberedThrough     
\ECRepeatTheorems

\EquationsNumberedThrough    

\usepackage{hyperref}

\begin{document}


\RUNAUTHOR{Guo, Hu, Xu and Zhang}

\RUNTITLE{Consistency and Computation of  Regularized MLEs for MHPs}

\TITLE{Consistency and Computation of  Regularized MLEs for Multivariate Hawkes Processes}


\ARTICLEAUTHORS{
\AUTHOR{Xin Guo}
\AFF{Department of Industrial Engineering and Operations Research, University of California, Berkeley, USA. 
\EMAIL{xinguo@berkeley.edu.}}
\AUTHOR{Anran Hu}
\AFF{Department of Industrial Engineering and Operations Research, University of California, Berkeley, USA. 
\EMAIL {anran\_hu@berkeley.edu.}}
\AUTHOR{Renyuan Xu}
\AFF{Department of Industrial Engineering and Operations Research, University of California, Berkeley, USA. 
\EMAIL {renyuanxu@berkeley.edu.}}
\AUTHOR{Junzi Zhang}
\AFF{Institute for Computational \& Mathematical Engineering, Stanford University, Stanford, USA. 
\EMAIL {junziz@stanford.edu}}
}

\ABSTRACT{
This paper proves the consistency property  for the regularized maximum likelihood estimators (MLEs) of multivariate Hawkes processes (MHPs). It also develops an alternating minimization type algorithm (AA-iPALM) to compute the MLEs with guaranteed global convergence to the set of stationary points. The performance of this AA-iPALM algorithm on both synthetic and real-world data shows that AA-iPALM consistently improves over iPALM and PALM. Moreover, AA-iPALM is able to  identify the causality structure in rental listings on craigslist  and herd behavior in the NASDAQ ITCH dataset.

}


\maketitle

%


\section{Introduction}

\paragraph{MHPs.} Multivariate Hawkes processes (MHPs) are counting/point  processes originally introduced in
\cite{HAWKES1971a,HAWKES1971b}  to model arrival patterns of earthquakes and the aftershocks triggered by
the earthquakes. There is a recent surge of interests in  MHPs to model social networks (\cite{BBH2012}, \cite{XFZ2016}, and \cite{YEHK2017}), earthquake events (\cite{OGATA1988a}, \cite{OGATA1998b}), and algorithmic tradings (\cite{HH2007},  \cite{MP2011}, \ 
\cite{EBK2012}, \cite{ZRA2014}, \cite{AJ2015},  \cite{BMM2015}, \cite{BJM2016},  
and \cite{GZZ2018}).

MHPs are characterized by their intensity processes, which represent the instantaneous likelihood of 
event arrivals. The intensity process has two components: the baseline intensity  and the triggering functions.
The former is similar to the intensity for a simple Poisson process, and the latter carry the mutual exciting property and measure how much the occurrence of one event ``triggers'' the arrival of others. 
The popularity of MHPs comes from this mutual exciting property: in social networks, triggering functions indicate how activities of one user  affect activities of other users; in financial markets, triggering functions are useful for capturing the herd behavior of trading activities.

\paragraph{Regularized MLEs.}

The standard  parameter estimation approach for MHPs is the maximum likelihood estimator (MLE). To encourage certain patterns and structures of the true parameter,  regularized MLEs have been proposed by adding a penalty term to the standard MLEs. The most common penalty terms include the $l_1$-norm in Lasso, the $l_2$-norm in Ridge, the
$l_{1,2}$-norm for group sparsity, and the nuclear norm for low-rank structure.
Regularized MLEs have been shown to be numerically efficient (\cite{ZZS2013_ADMM, ZZS2013_ode, BGM2015, XFZ2016}, and \cite{YEHK2017}). 

 In this paper, we will establish the consistency of regularized MLEs for MHPs. 
Our approach is inspired by  the deep work of \cite{OGATA1978},  and  extensively exploits the intriguing  connection between MHPs and the Branching processes. This connection first appeared in \cite{cluster_rep0} and was  ingeniously adapted by \cite{JHR2015} to study the higher order cumulants of MHPs. This powerful connection  enables us to avoid the abstract and technical assumptions used in \cite{OGATA1978}.

\paragraph{AA-iPALM algorithm.}

A well-recognized issue for learning MHPs is the computational challenge. Recently, \cite{NPEvsMLE} surveyed and  compared  MLEs with several  nonparametric approaches. Their study  confirmed the superb performance of 
MLEs for the one-dimensional Hawkes processes in terms of accuracy, and pointed out  the sensitivity of MLEs to initialization and hyper-parameter choices. Their survey highlighted the difficulty to deal with the nonconcave log-likelihood function for high dimensional Hawkes processes, and more importantly, the lack of  learning algorithms with guaranteed optimality in such scenarios.

 In this paper,  we propose an 
Anderson Accelerated inertial Proximal Alternating Linearized Minimization (AA-iPALM) algorithm to compute the MLEs of MHPs.
This algorithm exploits the block structure of MHPs, and combines iPALM (\cite{iPALM}) with a variant of Anderson acceleration technique recently proposed by \cite{AA1}.  We establish the global convergence of the AA-iPALM algorithm to the set of stationary points, and also obtain its iterative complexity. 

We test this AA-iPALM algorithm using both synthetic and real-world data, and compare its performance with both  iPALM and its non-inertial version PALM (proposed in \cite{PALM}). We show that AA-iPALM consistently improves over iPALM and PALM. Moreover, the AA-iPALM algorithm manages to identify some interesting causality structures in rental listings on craigslist  and herd behavior in the NASDAQ ITCH dataset.

\paragraph{Related works.}

 In addition to the well-known work of  \cite{OGATA1978} on the consistency of MLEs of general point processes,  \cite{PT1986} and \cite{GS2018} also studied the consistency of MLEs. The former focused 
 on  a special class of point processes with an almost surely bounded intensity process, and the latter assumed the Markovian property of the intensity process. 
Moreover, \cite{OZAKI1979} studied  the computation for MLEs of one-dimensional Hawkes process, and \cite{Clinet-Yoshida2017} established  the convergence of moments and hence asymptotic normality and consistency of MLEs
for exponential MHPs. However, none of these works considered the regularized MLEs for MHPs.

Apart from the MLE, there  is also the least-squares estimator (LSE) (see for example \cite{BGM2015}). Besides the parametric approach, the nonparametric approach has been reported in ~\cite{EM_mhp},~\cite{RASMUSSEN2013},~\cite{ZZS2013_ode},   ~\cite{BM2014},~\cite{LA2014}, ~\cite{DFASS2015},  ~\cite{kernel_parsim}, ~\cite{INAR}, ~\cite{XFZ2016}, ~\cite{EDD2017}, and ~\cite{YEHK2017}.

The iPALM algorithm (Algorithm \ref{alg:iPALM}) was proposed to solve non-convex optimization problems with block structures. It originated from  the Proximal Alternating Minimization algorithm
 (\cite{ABRS2010}), which was then developed into the Proximal Alternating Linearized Minimization algorithm (PALM) by replacing the exact minimization in the subproblems with a single gradient descent step  (\cite{PALM}). The iPALM algorithm  added additional inertia/momentum terms in each gradient descent step. 

 Anderson acceleration (AA) is an acceleration scheme for general fixed point problems (\cite{FangSecant, WalkerNi}). The main idea of AA is similar to limited memory quasi-Newton (QN) methods, and is also related to extrapolation methods. It has been applied to accelerate expectation maximization (EM) algorithms for computing MLEs of classical statistical models.

\section{Problem settings}\label{prob_set}


\paragraph{MHPs and intensity processes.}
A $K$-dimensional  MHP is a $K$-dimensional point process $\pmb{N}=(N_1, \cdots, N_K)$ defined on $\Omega:=\{\omega\;|\;\omega=\{t_j\;|\;j=0,\pm1,\pm2,\dots\}\mbox{ with no limit point}\}$. Here each $N_i(a,b]=N_i((a,b],\omega)$ is a one-dimensional point process counting the number of arrivals of the $i$-th type within the interval $(a,b]$. More precisely, $N_i(A)=N_i(A,\omega)$ is a counting measure defined for each bounded Borel set $A$ of $\mathbb{R}$, with $N_i(A,\omega)=\#(\omega\cap A)$. For notational simplicity, we also use the shorthand $N_i(t):=N_i(0,t]$ in some cases.

Throughout the paper, we define the filtration $\mathcal{F}_t$ to be the $\sigma$-field generated by $\{\pmb{N}(-\infty,s]\;|\; s\leq t\}$. We say a process $\xi=\{\xi(t,\omega)\;|\;t\geq 0\}$ is \textit{adapted} if $\xi(t,\omega)$ is $\mathcal{F}_{t^-}$-measurable for every fixed $t\geq 0$. It is said to be {\it predictable} if $\xi$, considered as a mapping from $\mathbb{R}_+\times \Omega$ to $\mathbb{R}$, is measurable with respect to the $\sigma$-algebra generated by all left continuous adapted processes.  

An MHP is characterized by its intensity process. 
There are two standard definitions of intensity processes for MHPs in the literature.  (See \cite{LINIGER2009}.)

The first definition is from the original work of \cite{HAWKES1971a, HAWKES1971b}, where the intensity process $\pmb{\lambda} = (\lambda_1, \cdots, \lambda_K)$  is given as the following. For each $i=1,\cdots,K$,  
\begin{equation}\label{trueint}
\lambda_{i}(t;\pmb{\theta})=\mu_i+\sum_{j=1}^K\int_{-\infty}^{t-} g_{ij}(t-s;\bigeta)N_j(ds). 
\end{equation}
Here $\btheta=(\bmu,\bigeta)=(\mu_1,\cdots,\mu_K,\eta_1,\dots,\eta_D)\in\Theta\subseteq\mathbb{R}^{K+D}$ is the underlying parameter with $\mu_i$  the baseline intensity of $N_i$,  and $g_{ij}:\mathbb{R}_+\times \mathbb{R}^D\rightarrow \mathbb{R}$ the triggering function that captures the mutual excitation or the causality between $N_i$ and $N_j$ for $i\ne j$ and the self excitation for $i=j$.

The second definition appears in most of the engineering literature, where the intensity process $\hat{\blambda} = (\hat{\lambda}_1, \cdots, \hat{\lambda}_K)$ takes the following form 
\begin{equation}\label{int}
\hat{\lambda}_{i}(t;\btheta)=\mu_i+\sum_{j=1}^K\int_{0}^{t-}g_{ij}(t-s;\bigeta)N_j(ds).
\end{equation}

Note that the first definition is infeasible for parameter estimations, as data from an infinite time interval is inaccessible. However, this intensity process has nice mathematical properties such as stationarity
under mild conditions on $g_{ij}$ and $\mu_i$. 
The second definition is more natural for parameter estimations. However, this intensity process, unless assumed
a constant,  is nonstationary. In fact, lack of the stationarity property is the main difficulty in establishing the consistency for MLEs.

\paragraph{Corresponding MLEs.}
 The log-likelihood functions for the intensity processes (\ref{trueint}) and (\ref{int})  on the interval $[0,T]$ take the following forms, respectively,
\begin{equation}\label{exact_l}
L_T(\btheta) := \sum_{i=1}^K\left[-\int_0^T \lambda_i(t;\btheta) dt + \int_0^T\log(\lambda_i(t; \btheta)) N_i(dt)\right],
\end{equation}
\begin{equation}\label{exact_ll}
\hat{L}_T(\btheta) := \sum_{i=1}^K\left[-\int_0^T \hat{\lambda}_i(t;\btheta) dt + \int_0^T\log(\hat{\lambda}_i(t; \btheta)) N_i(dt)\right].
\end{equation}
Their respective MLEs are 
\begin{equation}\label{mle}
\btheta_T:=\argmax\nolimits_{\btheta\in\Theta}L_T(\btheta),\quad\hat{\btheta}_T:=\argmax\nolimits_{\btheta\in\Theta}\hat{L}_T(\btheta).
\end{equation}

The regularized MLE for $\hat{L}_T(\btheta)$ is defined by adding a penalty term $P(\btheta)$, i.e., 
\begin{equation}\label{mle_reg}
\hat{\btheta}_T^{reg}:=\argmax\nolimits_{\btheta\in\Theta}\hat{L}_T(\btheta)-P(\btheta)=\argmax\nolimits_{\btheta\in\Theta}\hat{L}_T^{reg}(\btheta).
\end{equation}

For notational simplicity, we will sometimes use the shorthand $\lambda_i(t):=\lambda_i(t;\btheta^\star)$ to denote the intensity process corresponding to the underlying true MHP $\pmb{N}$. 

\section{Consistency of MLEs}\label{asymptotic}

In this section, we will establish the consistency for $\hat{\btheta}_T^{reg}$. That is,  it converges in probability to the (unknown) true parameter, denoted as {${\btheta}^{\star}=(\bmu^\star,\pmb{\eta}^\star)\in\Theta$}. This consistency also implies the convergence of the vanilla MLE $\hat{\btheta}_T$.

\begin{theorem}[Consistency] \label{parameter_conv}
Under Assumptions \ref{assumption1}, \ref{assumption2}, \ref{assumption4} and \ref{assumption3} (specified below), and assume that $P$ is continuous in $\btheta$, the regularized MLE $\hat{\btheta}_T^{reg}$  converges to ${\btheta}^{\star}$ in probability as $T \rightarrow \infty$.
\end{theorem}

\subsection{Preliminaries}

The first critical piece is the Branching process representation of an MHP.

\paragraph{Branching process representation for MHPs.}
 As \cite{cluster_rep0}  pointed out, other than its intensity process, a Hawkes process can be equivalently  defined as a Poisson cluster process with a certain branching structure. More precisely, a $K$-dimensional MHP with positive baseline intensity $\bmu$ and nonnegative integrable triggering functions $g_{ij}$ can be equivalently constructed as the following process. (See \cite{JHR2015}, and also \cite{RASMUSSEN2013} for the one-dimensional cases.)
 \begin{itemize}
\item For $k=1,2,\dots,K$, initialize an instance $I_k$ of a homogeneous Poisson process with rate $\mu_k$, with its elements called immigrants of type $k$;
\item Immigrants of all types generate independent clusters. More specifically, for each $k=1,2,\dots,K$, each immigrant $x\in I_k$ generates a cluster $C_x^k$ with the following branching structure: 
\begin{itemize}
\item Generation $0$ consists of the immigrant $x$;
\item Recursively, given generations $0,1,\dots,n$, for all $i,j=1,\dots,K$, each point $s$ of type $j$ in generation $n$ generates its offsprings of type $i$ as an instance of an inhomogeneous Poisson process with rate $\lambda_{ij}(t):=g_{ij}(t-s)$. All these offsprings then constitute generation $n+1$.
\end{itemize}
\item The point process is then defined to be the union of all clusters.
\end{itemize}

Under the original notation of MHPs, the number of type $i$ points (immigrants and offsprings) within time interval $(a,b]$ is exactly $N_i(a,b]$.

Next, we recall some basic results in analysis.

The first proposition, stated informally in \cite{OGATA1978},  translates a stochastic Lebesgue-Stieltjes integral to a Lebesgue integral.

\begin{proposition}\label{dN2dt_ogata}
Suppose that $\{\xi(t);~t\geq 0\}$ is a finite predictable process such that 
$\int_0^T\mathbb{E}\left[|\xi(t)|\lambda_i(t;\btheta^\star)\right]dt<\infty$ for any $T\geq 0$, then
\[
\mathbb{E}\left[\int_0^T\xi(t)N_i(dt)\right]=\mathbb{E}\left[\int_0^T\xi(t)\lambda_i(t;\btheta^{\star})dt\right].
\]
\end{proposition}

The next proposition controls the mean-square difference between the Lebesgue-Stieltjes integral and the Lebesgue integral. {(See \cite{Protter}, Theorem 20 in Section II.5, Theorem 29 and  Corollary 3 in Section II.6)}.

\begin{proposition}\label{dN2dt}
Suppose that $\{\xi(t);~t\geq 0\}$ is a finite predictable process, with $\mathbb{E}\left[\int_0^T\xi(t)^2N_i(dt)\right]<\infty$ for any $T\geq 0$, $i=1,\dots,K$. Define $M_i(t):=N_i(t)-\int_0^t\lambda_i(s;\btheta^\star)ds$ and $X_i(t):=\int_0^t\xi(s)M_i(ds)$. Then 
\[
\mathbb{E}[X_i(t)^2]=\mathbb{E}\left[\int_0^t\xi(s)^2N_i(ds)\right].
\]
\end{proposition}

The following result connects the convergence of a sequence with the convergence of its C\'esaro-sum sequence.

\begin{proposition}\label{cesaro}
Suppose that $f(t):\mathbb{R}_+\rightarrow\mathbb{R}_+$ is Lebesgue measurable, and $f(t)\leq C$ for all $t\geq 0$ and some $C>0$. If $f(t)\rightarrow0$ as $t\rightarrow\infty$, then $\frac{1}{T}\int_0^Tf(t)dt\rightarrow 0$ as $T\rightarrow\infty$.
\end{proposition}

\subsection{General assumptions and key lemmas}\label{gen_assump}
Next, we   introduce the  assumptions that will be used in this section for the consistency of MLEs.
For simplicity, we denote $\Theta_{\pmb{\eta}}:=\{\pmb{\eta}\;|\;\exists\,\bmu, (\bmu,\pmb{\eta})\in\Theta\}$. 

The first set  of technical conditions ensures the regularity of the intensity process and the log-likelihood functions. 
In particular, the compactness of $\Theta$ defines  a rough range of the true parameters, and the bounds on $g_{ij}$ ensure the well-definedness for  the consistency  of MLEs.

\begin{assumption}[Regularity]\label{assumption1}
$\Theta$ is nonempty and compact. For $i=1,\dots,K$, $\mu_i\geq \underline{\mu}$ for some $\underline{\mu}>0$. For $i, j=1,\cdots,K$, $g_{ij}$ is bounded, nonnegative, {left continuous and integrable over $[0,\infty)$ with respect to $t$ for any $\bigeta\in \Theta_{\bigeta}$}. Moreover, there is an open set $\tilde{\Theta}\supseteq \Theta$ such that $g_{ij}(t;\pmb{\eta})$ is differentiable with respect to $\pmb{\eta}$ in $\tilde{\Theta}$. 
\end{assumption}

The second assumption on stationarity is standard  in the  existing literature.
It ensures that the MHP and its underlying intensity process $\pmb{\lambda}(t;\btheta^\star)$ are stationary and ergodic.
In particular, $\bar{{\lambda}}_i:=\mathbb{E}[{\lambda}_i(t;\pmb{\theta^\star})]$ is  independent of $t$. 
(See Theorem 7 in \cite{Bremaud-Massoulie}.)

\begin{assumption}[Stationarity]\label{assumption2}
The spectral radius of matrix $\pmb{G}:=[G_{ij}]_{K \times K}$ is smaller than $1$, where $G_{ij}:=\int_0^{\infty}g_{ij}(t;{\pmb{\eta}}^{\star})dt$.
\end{assumption}
Given Assumptions \ref{assumption1} and \ref{assumption2}, 
the aforementioned connection between the Branching process and an MHP leads to an explicit formula for cumulants of arbitrary orders, which is crucial for the  following lemma concerning moment bounds on the number of arrivals within given intervals.  This lemma is critical for all the following lemmas.
\begin{lemma}[Main lemma]\label{moment}
Given Assumptions \ref{assumption1} and \ref{assumption2},
$\exists ~\bar{C}>0$ such that \[\max_{i=1,\dots,K}\mathbb{E}\left[\left|{N_i}(t,t+h]-\bar{{\lambda}}_ih\right|^{4}\right]\leq \bar{C}h^{3},\] 
$\forall~ t\in\mathbb{R},~ h>0$, where $\bar{{\lambda}}_i:=\mathbb{E}[{\lambda}_i(t;\pmb{\theta^\star})]$. 
\end{lemma}

The next  assumption  is essential for distinguishing the true parameter from the others. 
\begin{assumption}[Identifiability]\label{assumption4}
For any $\pmb{\eta}\neq \pmb{\eta}'\in \Theta_{\pmb{\eta}}$, there exists a set of $t\in\mathbb{R}_+$ with nonzero measure such that $\pmb{G}(t;\pmb{\eta})\neq \pmb{G}(t;\pmb{\eta}')$. Here $\pmb{G}(t;\pmb{\eta}):=[g_{ij}(t;\pmb{\eta})]_{K\times K}$.
\end{assumption}

\begin{lemma}\label{identifiability}
Under Assumptions \ref{assumption1}, \ref{assumption2} and \ref{assumption4},  $\blambda(0;\btheta)=\blambda(0;\btheta')$ a.s. if and only if $\btheta=\btheta'$.
\end{lemma}

The summability condition below  requires  that the triggering functions  decay sufficiently fast. 

\begin{assumption}[Summability]\label{assumption3}
For any $i,j=1,\cdots,K$, $d=1,\dots,D$, $g_{ij}(t;\pmb{\eta})$ and its partial derivatives $\partial_{\eta_d}g_{ij}$ are all uniformly summable. 
In addition
\[
\lim_{t\rightarrow\infty}\sum_{k=1}^{\infty}(t_k-t_{k-1})\sup\nolimits_{t'\in[t+t_{k-1},t+t_k],\pmb{\eta}'\in \Theta_{\pmb{\eta}}}g_{ij}(t';\pmb{\eta}')=0,
\]
where the sequence $\{t_k\}_{k=0}^{\infty}$ satisfies $t_0=0$ and $\sum_{k=1}^{\infty}(t_k-t_{k-1})^{-1}<\infty$. Here 
a function $h(t;\pmb{\eta})$, $h:\mathbb{R}\times\Theta_{\pmb{\eta}}\rightarrow\mathbb{R}$ is uniformly summable if there exists a strictly increasing sequence $\{t_k\}_{k=0}^{\infty}$ such that for all
$t\geq 0$,
\[
\begin{split}
\sum_{k=1}^{\infty}(t_k-t_{k-1})\sup\nolimits_{t'\in[t+t_{k-1},t+t_k],\pmb{\eta}'\in\Theta_{\pmb{\eta}}}|h(t';\pmb{\eta}')|<E,
\end{split}
\]
for some constant $E>0$, 
where the sequence $\{t_k\}_{k=0}^{\infty}$ satisfies $t_0=0$ and $\sum_{k=1}^{\infty}(t_k-t_{k-1})^{-1}<\infty$.

\end{assumption}

 It is easy to verify  Assumption 
\ref{assumption3} for standard triggering functions, including the exponential, the Rayleigh, and the power-law types. For instance, 
let $t_0=0$, $t_k-t_{k-1}=k^{1+\epsilon'}$ ($k\geq 1$) for some $\epsilon'>0$, then $\sum_{k=1}^{\infty}(t_k-t_{k-1})^{-1}=\sum_{k=1}^{\infty}k^{-1-\epsilon'}<\infty$, then Assumption \ref{assumption3} holds for all triggering functions with a polynomial-exponential decay of the form $t^me^{-\beta t}$ with $m,\beta>0$. 
Moreover, as the summability property is closed under summations, Assumption \ref{assumption3} is satisfied for most triggering functions in the literature. 

Combining Assumptions \ref{assumption1}, \ref{assumption2}, and
\ref{assumption3}, we can derive some useful bounds and  the asymptotic difference between the intensity processes ${\blambda}$ and  $\hat{\blambda}$. 

 \begin{lemma}\label{uniform_increment}
Under Assumptions \ref{assumption1}, \ref{assumption2} and \ref{assumption3}, for any $t\geq 0$, define for $i=1,\dots,K$ that
\[C_i^{(t)}:=\sup\nolimits_{k\geq 1}N_i(t-t_k,t-t_{k-1}]/(t_k-t_{k-1}).
\]
Then $\mathbb{E}[|C_i^{(t)}|^{3+\alpha}]\leq C$ for any $\alpha\in[0,1)$. Here $C>0$ is a positive constant independent of $t$, and $\{t_k\}_{k=0}^{\infty}$ is the sequence from Assumption \ref{assumption3}.
\end{lemma}

Note that Lemma \ref{uniform_increment} follows directly from Lemma \ref{moment},  and is critical for Lemmas \ref{newlemma_consist}, \ref{approx}, and  \ref{cont_lambda}.

\begin{lemma}\label{newlemma_consist}
Under Assumptions \ref{assumption1}, \ref{assumption2} and \ref{assumption3}, there exist two random variables $\Lambda_0,~\Lambda_1$ with finite $(3+\alpha)$-th moments for any $\alpha\in[0,1)$, such that for any $1\leq i\leq K$,
\[
\sup\nolimits_{\pmb{\theta}'\in \Theta}\lambda_i(0;\pmb{\theta}')\leq \Lambda_0,\quad \sup\nolimits_{\pmb{\theta}'\in \Theta}|\log \lambda_i(0;\pmb{\theta}')|\leq \Lambda_1.
\]

\end{lemma}

\begin{lemma}\label{approx}
 Under Assumptions \ref{assumption1}, \ref{assumption2} and \ref{assumption3}, for all $i=1,\dots,K$,
  \[
  \begin{split}
 &  \text{$\mathbb{E}\left[\sup\nolimits_{\btheta\in \Theta} \left|\lambda_i(t;\btheta)-\hat{\lambda}_i(t;\btheta)\right|^2\right]\leq F$ for all $t\geq 0$ and some constant $F>0$,}\\
 & \text{$\mathbb{E}\left[\sup\nolimits_{\btheta\in \Theta} \left|\lambda_i(t;\btheta)-\hat{\lambda}_i(t;\btheta)\right|^2\right]\rightarrow 0$ as $t\rightarrow\infty$}.
 \end{split}
\]
\end{lemma}

\begin{lemma}\label{cont_lambda}
Under Assumptions \ref{assumption1}, \ref{assumption2} and \ref{assumption3},
 $\blambda(t;\btheta)$ is a.s. continuous in $\btheta$ for any fixed $t\geq 0$.
 \end{lemma}

Furthermore,  the stationarity and ergodicity of the underlying MHP, together with Proposition \ref{dN2dt_ogata} and Proposition \ref{dN2dt}, lead to the following lemma.
\begin{lemma}\label{ergodic_xi}
Suppose that $\xi=\{\xi(t);t\geq 0\}$ is a stationary stochastic process, {and suppose that it is adapted and  a.s. left continuous in $t\geq 0$}, with $\mathbb{E}[|\xi(0)|^2]<\infty$. 
Then the following limits hold:
\begin{equation}\label{ergodic_xi_1} 
\frac{1}{T}\int_0^T\xi(t)dt\overset{\mathbb{P}}{\longrightarrow}\mathbb{E}[\xi(0)],
\end{equation}

\begin{equation}\label{ergodic_xi_2}
\frac{1}{T}\int_0^T\dfrac{\xi(t)}{\lambda_i(t;\btheta^\star)}N_i(dt)\overset{\mathbb{P}}{\longrightarrow}\mathbb{E}[\xi(0)].
\end{equation}
\end{lemma}

In contrast to  Lemma 2 in \cite{OGATA1978} which implicitly assumes the almost sure convergence of $\frac{1}{T}\int_{\lfloor T\rfloor }^T\xi(t)dt$ to $0$ to establish the almost sure convergence of (\ref{ergodic_xi_1}) and (\ref{ergodic_xi_2}),  the above  property of convergence in probability suffices for the consistency of MLEs.

Also, one can establish the predictability, stationarity and certain bounds and approximation errors between the log-likelihood functions. 
 \begin{lemma}\label{stat_moment_inflambda}
 Under Assumptions \ref{assumption1}, \ref{assumption2} and \ref{assumption3},
for any $\btheta\in \Theta$, $\lambda_i(t;\btheta)$ is stationary, adapted, and a.s. left continuous in $t\geq 0$. Moreover, for any subset $U\subseteq \Theta$, the following stochastic processes
 \[
 \left\{
 \begin{array}{l}
 \xi_{U,i}^{(1)}(t):=\inf_{\btheta\in U}\lambda_i(t;\btheta)-\lambda_i(t;\btheta^\star),\\
 \xi_{U,i}^{(2)}(t):=\lambda_i(t;\btheta^\star)\left(\log\left(\lambda_i(t;\btheta^\star)\right)-\log\left(\sup_{\btheta\in U}\lambda_i(t;\btheta)\right)\right),\\ 
 \xi_{U,i}^{(3)}(t):=\sup_{\btheta\in U}\left(\lambda_i(t;\btheta)-\hat{\lambda}_i(t;\btheta)\right)
 \end{array}
 \right.\qquad i=1,\dots,K,
 \]
 are also adapted and a.s. left continuous in $t\geq 0$. In addition, $\xi_{U,i}^{(1)}(t)$ and $\xi_{U,i}^{(2)}(t)$ are stationary, and  for any $T\geq 0$ and $i=1,\dots,K$, 
 \[
 \begin{array}{l}
\mathbb{E}\left[\left|\xi_{U,i}^{(1)}(0)\right|^2\right]<\infty,
\quad \mathbb{E}\left[\left|\xi_{U,i}^{(2)}(0)\right|^2\right]<\infty,\quad \int_0^T\mathbb{E}\left[\left|\xi_{U,i}^{(3)}(t)\right|\lambda_i(t;\btheta^\star)\right]dt<\infty.
 \end{array}
 \]
 \end{lemma}
 It is worth noticing that
if $\{\xi(t);~t\geq 0\}$ is an adapted stochastic process defined on the filtered probability space of the (true) MHP, and if the sample paths of $\xi(t)$ are a.s. left continuous on $(0,\infty)$, then $\xi(t)$ is predictable. Thus the aforementioned stochastic processes $\lambda_i(t;\btheta)$, $\xi_{U,i}^{(l)}(t)$, $l=1,2,3$ are all predictable.

 By plugging $\xi_{U,i}^{(1)}(t)$ in Eqn. (\ref{ergodic_xi_1}) and $\xi_{U,i}^{(2)}(t)$ in Eqn. (\ref{ergodic_xi_2}), respectively, 
Lemma \ref{ergodic_xi} and Lemma \ref{stat_moment_inflambda} imply the following ergodicity property.
\begin{lemma}[Ergodicity]\label{ergodic}
Under Assumptions \ref{assumption1}, \ref{assumption2} and \ref{assumption3}, for all $1\leq i\leq K $ and any subset $U\subseteq \Theta$, the following limits hold: 
\begin{equation}
 \frac{1}{T} \int_{0}^T \left(\inf_{\btheta\in U}\lambda_i(t;\btheta)-\lambda_i(t;{\btheta}^{\star})\right) dt  \overset{\mathbb{P}}{\longrightarrow} \mathbb{E}\left[\inf_{\btheta\in U}\lambda_i(0;\btheta)-\lambda_i(0;{\btheta}^{\star})\right],
\end{equation}

\begin{equation}
  \frac{1}{T} \int_0^T \log\left(\frac{\lambda_i(t;{\btheta}^{\star})}{\sup_{\btheta\in U}\lambda_i(t;\btheta)}\right) N_i(dt)
\overset{\mathbb{P}}{\longrightarrow}
\mathbb{E}\left[\lambda_i(0;\pmb{\theta^\star})\log\left(\frac{\lambda_i(0;{\btheta}^{\star})}{\sup_{\btheta\in U}\lambda_i(0;\btheta)}\right)\right].
\end{equation}
\end{lemma}

We will also need the following result related to a Kullback-Leibler's divergence type quantity defined from the intensity processes $\blambda(t;\btheta)$.
\begin{lemma}\label{KL}
Under Assumptions \ref{assumption1}, \ref{assumption2}, \ref{assumption4} and \ref{assumption3}, 
\begin{equation}
D({\btheta}^{\star};\btheta):=\sum_{i=1}^K\mathbb{E}\left[\lambda_i(0;{\btheta}^{\star})\left\{\frac{\lambda_i(0;\btheta)}{\lambda_i(0;{\btheta}^{\star})}-1+\log\left(\frac{\lambda_i(0;\btheta^{\star})}{\lambda_i(0;{\btheta})}\right)\right\}\right]\geq 0.
\end{equation}
The equality holds if and only if\, $\btheta=\btheta^\star$. In addition, $D(\btheta^\star;\btheta)$ is continuous in $\btheta$ at $\btheta^\star$. 
\end{lemma}
The proof is a combination of Lemma \ref{cont_lambda} and Lemma \ref{identifiability}, with the classical dominated convergence theorem due to Lemma \ref{newlemma_consist}.

\begin{remark}\label{compare_ogata}
In essence, the Branching process representation of MHPs enables us to replace all the abstract assumptions in \cite{OGATA1978} with Assumptions \ref{assumption1},
\ref{assumption2}, \ref{assumption4}, and \ref{assumption3}. The consequent  Lemmas \ref{KL},  \ref{identifiability}, \ref{cont_lambda}, \ref{newlemma_consist}, and \ref{approx} in this paper replace the  continuity assumption of $D(\btheta^\star;\btheta)$, Assumption B3, 
Assumption B4, Assumption B5, and  Assumption C1(i) in \cite{OGATA1978}, respectively.
\end{remark}

\subsection{Proof of Theorem \ref{parameter_conv}}
Our proof adapts the proof scheme in \cite{OGATA1978}  to our problem setting for regularized MLEs, with all details.
It can be divided into two steps. The first step is to show   the consistency of $\btheta_T$. That is,  for any neighborhood $U_0$ of $\btheta^\star$,  $\exists$ $\epsilon>0$, such that as $T\rightarrow\infty$,
 \begin{eqnarray}\label{ine_2_main}
\mathbb{P}\left(\sup\nolimits_{\btheta \in U_0} L_T(\btheta) \geq \sup\nolimits_{\btheta \in \Theta \setminus U_0} L_T(\btheta) + \epsilon T\right)\rightarrow 1.
\end{eqnarray}
This follows from  Lemma \ref{ergodic} and  Lemma \ref{KL}.

The second step is to utilize the closeness of $\pmb{\lambda}$ and $\hat{\pmb{\lambda}}$ to show that
\begin{eqnarray}\label{ine_2_main'}
\mathbb{P}\left(\sup\nolimits_{\btheta \in U_0} \hat{L}_T^{reg}(\btheta) \geq \sup\nolimits_{\btheta \in \Theta \setminus U_0} \hat{L}_T^{reg}(\btheta) + \epsilon T/4 \right)\rightarrow 1.
\end{eqnarray}
This can be established by Lemmas \ref{newlemma_consist},  \ref{approx},  \ref{stat_moment_inflambda}, and Propositions \ref{dN2dt_ogata}  and \ref{cesaro}. 
The details are below.

\paragraph{Step 1: Consistency for $L_T$ and $\btheta_T$.}
Let $U$ be any neighborhood of $\btheta$. When $U$ shrinks to $\{\btheta\}$, by Lemma \ref{newlemma_consist} and Lemma \ref{cont_lambda}, we can apply the Lebesgue dominated convergence theorem to obtain  
\[
\mathbb{E}  \left[\sum_{i=1}^K\inf_{\btheta' \in U}\lambda_i(0;\btheta')\right] \rightarrow \mathbb{E} \left[\sum_{i=1}^K \lambda_i(0; \btheta)\right],
\]
and\[
\mathbb{E}  \left[\sum_{i=1}^K \lambda_i(0;{\btheta}^{\star}) \log \left(\frac{\lambda_i(0;{\btheta}^{\star})}{\sup_{\btheta^{'} \in U} \lambda_i(0;\btheta^{'})} \right) \right] \rightarrow \mathbb{E}  \left[\sum_{i=1}^K \lambda_i(0;{\btheta}^{\star}) \log \left(\frac{\lambda_i(0;{\btheta}^{\star})}{\lambda_i(0;\btheta)} \right)\right].
\]
Here we use the fact that for any $U\subseteq\Theta$,
\[
\left|\sum_{i=1}^K\inf_{\btheta' \in U}\lambda_i(0;\btheta')\right|\leq K\Lambda_0,\quad \left|\sum_{i=1}^K \lambda_i(0;{\btheta}^{\star}) \log \left(\frac{\lambda_i(0;{\btheta}^{\star})}{\sup_{\btheta^{'} \in U} \lambda_i(0;\btheta^{'})} \right)\right| \leq 2K\Lambda_0\Lambda_1, 
\]
$
\mathbb{E}[\Lambda_0]<\infty$, and $\mathbb{E}[2K\Lambda_0\Lambda_1]\leq 2K\sqrt{\mathbb{E}[\Lambda_0^2]\mathbb{E}[\Lambda_1^2]}<\infty$
implied from Lemma \ref{newlemma_consist}.

Let $U_0$ be an open neighborhood of ${\btheta}^{\star}$. Then by Lemma \ref{KL}, there exists $\epsilon>0$ such that $D({\btheta}^{\star};\btheta) \geq 3 \epsilon$  for any $\btheta \in \Theta \setminus U_0$. Now for any $\btheta \in \Theta \setminus U_0$, one can choose a sufficiently small  open neighborhood $U_{\btheta}$ of $\btheta$ such that
\begin{eqnarray*}
\mathbb{E} \left[\sum_{i=1}^K \inf_{\btheta'\in U_{\btheta}} \lambda_i(0;\btheta') - \sum_{i=1}^K\lambda_i(0;{\btheta}^{\star}) + \sum_{i=1}^K \lambda_i(0;{\btheta}^{\star}) \log \left(\frac{\lambda_i(0;{\btheta}^{\star})}{\sup_{\btheta' \in U_{\btheta}} \lambda_i(0;\btheta')}\right)\right]\geq  D({\btheta}^{\star};\btheta)-\epsilon \geq 2\epsilon.
\end{eqnarray*}

By the finite covering theorem, one can select a finite number of $\btheta_s\in \Theta\setminus U_0$, $1 \leq s \leq N$ such that the union of the sets $U_{s}=U_{\btheta_s}$ covers $\Theta \setminus U_0$. By Lemma \ref{ergodic}, for $s = 1,2,\cdots,N$, 

\begin{align}
\lim_{T\rightarrow\infty}\mathbb{P}&\left(\frac{1}{T}\int_0^T \left(\sum_{i=1}^K \inf_{\btheta \in U_s}\lambda_i(t;\btheta)- \lambda_i(t;{\btheta}^{\star}) \right)dt + \frac{1}{T} \sum_{i=1}^K  \int_0^T \log \left(\frac{\lambda_i(t;{\btheta}^{\star})}{\sup_{\btheta \in U_s}\lambda_i(t;\btheta)}\right)N_i(dt) \nonumber \right.\\
 &\left.\geq\mathbb{E} \left[\sum_{i=1}^K \inf_{\btheta'\in U_{s}} \lambda_i(0;\btheta') - \sum_{i=1}^K\lambda_i(0;{\btheta}^{\star}) + \sum_{i=1}^K \lambda_i(0;{\btheta}^{\star}) \log \left(\frac{\lambda_i(0;{\btheta}^{\star})}{\sup_{\btheta' \in U_{s}} \lambda_i(0;\btheta')}\right)
 \right]-\epsilon
\right)=1,
\end{align}
which implies 
\begin{equation}\label{consist_L}
\lim_{T\rightarrow\infty}\mathbb{P}\left(\frac{1}{T}\int_0^T \left(\sum_{i=1}^K \inf_{\btheta \in U_s}\lambda_i(t;\btheta)- \lambda_i(t;{\btheta}^{\star}) \right)dt + \frac{1}{T} \sum_{i=1}^K  \int_0^T \log \left(\frac{\lambda_i(t;{\btheta}^{\star})}{\sup_{\btheta \in U_s}\lambda_i(t;\btheta)}\right)N_i(dt)\geq\epsilon \right)=1.
\end{equation}

Since 
\begin{equation}\label{ine_1}
\begin{split}
\frac{1}{T} L_T ({\btheta}^{\star}) &- \sup_{\btheta \in U_s} \frac{1}{T} L_T(\btheta)\\
&\geq  \frac{1}{T}\int_0^T \left(\sum_{i=1}^K \inf_{\btheta \in U_s}\lambda_i(t;\btheta)- \lambda_i(t;{\btheta}^{\star}) \right)dt + \frac{1}{T} \sum_{i=1}^K  \int_0^T \log \left(\frac{\lambda_i(t;{\btheta}^{\star})}{\sup_{\btheta \in U_s}\lambda_i(t;\btheta)}\right)N_i(dt),
\end{split}
\end{equation}
 $\lim_{T\rightarrow \infty}\mathbb{P}\left(\dfrac{1}{T} L_T ({\btheta}^{\star}) - \sup_{\btheta \in U_s}\frac{1}{T} L_T(\btheta)\geq \epsilon\right)=1$. Moreover, since $\btheta^\star\in U_0$, $N$ is finite and the union of $U_{\btheta_s}$ covers $\Theta\setminus U_0$, 
\begin{eqnarray}\label{ine_2}
\lim_{T\rightarrow\infty}\mathbb{P}\left(\sup_{\btheta \in U_0} L_T(\btheta) \geq \sup_{\btheta \in \Theta \setminus U_0} L_T(\btheta) + \epsilon T\right)=1.
\end{eqnarray} 
This implies $\lim_{T\rightarrow\infty}\mathbb{P}\left(\btheta_T\in U_0\right)=1$ and the consistency of $\btheta_T$ is then established by the arbitrariness of $U_0$.

\paragraph{Step 2. Consistency for $\hat{L}_T^{reg}$ and $\hat{\btheta}_T^{reg}$.} 
As in the first step, to establish the consistency of $\hat{\btheta}_T^{reg}$ is to show 
\begin{equation}\label{ine_approx_18}
\lim_{T\rightarrow \infty}\mathbb{P}\left(\dfrac{1}{T} \hat{L}_T^{reg} ({\btheta}^{\star}) - \sup_{\btheta \in U_s}\frac{1}{T} \hat{L}_T^{reg}(\btheta)\geq \epsilon/4\right)=1.
\end{equation}
Since $\Theta$ is compact and $P$ is continuous, we have $\sup_{\btheta\in \Theta}|P(\btheta)|<\infty$. Hence it suffices to prove that 
\begin{equation}\label{ine_approx_18_ver2}
\lim_{T\rightarrow \infty}\mathbb{P}\left(\dfrac{1}{T} \hat{L}_T ({\btheta}^{\star}) - \dfrac{1}{T}\sup_{\btheta \in U_s} \hat{L}_T(\btheta)-\epsilon/4\geq \dfrac{1}{T}P(\btheta^\star)+\dfrac{1}{T}\sup_{\btheta\in U_s}|P(\btheta)|\right)=1.
\end{equation}
By noticing that $\lim_{T\rightarrow\infty}\left|\frac{1}{T}P(\btheta^\star)+\frac{1}{T}\sup_{\btheta\in U_s}|P(\btheta)|\right|\leq \lim_{T\rightarrow\infty}\frac{2}{T}\sup_{\btheta\in\Theta}|P(\btheta)|=0$, we see that limit (\ref{ine_approx_18_ver2}) can be further reduced to the following:
\begin{equation}\label{ine_approx}
\lim_{T\rightarrow \infty}\mathbb{P}\left(\dfrac{1}{T} \hat{L}_T ({\btheta}^{\star}) - \sup_{\btheta \in U_s}\frac{1}{T} \hat{L}_T(\btheta)\geq \epsilon/2\right)=1.
\end{equation}

To prove this, notice that
\begin{equation}\label{ine_1'}
\begin{split}
\frac{1}{T} \hat{L}_T ({\btheta}^{\star}) &- \sup_{\btheta \in U_s} \frac{1}{T} \hat{L}_T(\btheta)\\
&\geq  \frac{1}{T}\int_0^T \left( \sum_{i=1}^K \inf_{\btheta \in U_s}\hat{\lambda}_i(t;\btheta)- \hat{\lambda}_i(t;{\btheta}^{\star}) \right)dt + \frac{1}{T} \sum_{i=1}^K  \int_0^T \log \left(\frac{\hat{\lambda}_i(t;{\btheta}^{\star})}{\sup_{\btheta \in U_s}\hat{\lambda}_i(t;\btheta)}\right)N_i(dt).
\end{split}
\end{equation}
By (\ref{consist_L}), we only need to show 
\begin{align}
\lim_{T\rightarrow\infty}&\mathbb{P}\left(\frac{1}{T}\int_0^T \left( \sum_{i=1}^K \inf_{\btheta \in U_s}\hat{\lambda}_i(t;\btheta)- \hat{\lambda}_i(t;{\btheta}^{\star}) \right)dt + \frac{1}{T} \sum_{i=1}^K  \int_0^T \log \left(\frac{\hat{\lambda}_i(t;{\btheta}^{\star})}{\sup_{\btheta \in U_s}\hat{\lambda}_i(t;\btheta)}\right)N_i(dt) \right. \nonumber \\
&\geq\left.\frac{1}{T}\int_0^T \left(\sum_{i=1}^K \inf_{\btheta \in U_s}\lambda_i(t;\btheta)- \lambda_i(t;{\btheta}^{\star}) \right)dt + \frac{1}{T} \sum_{i=1}^K  \int_0^T \log \left(\frac{\lambda_i(t;{\btheta}^{\star})}{\sup_{\btheta \in U_s}\lambda_i(t;\btheta)}\right)N_i(dt)-\epsilon/2\right)=1.
\end{align}

Recall that $\hat{\lambda}_i(t;\btheta)\leq\lambda_i(t;\btheta)$ for all $i=1,\dots,K$ and for any $t$ and $\btheta\in\Theta$, it suffices to show that the following limits hold as $T\rightarrow\infty$: 
\begin{align}
&\mathbb{E}\left[\frac{1}{T}\int_0^T \left[ \inf_{\btheta \in U_s}{\lambda}_i(t;\btheta)-\inf_{\btheta \in U_s}\hat{\lambda}_i(t;\btheta)\right]dt\right]\rightarrow 0\label{1},\\
&\mathbb{E}\left[\frac{1}{T}\int_0^T \left[ {\lambda}_i(t;\btheta^\star)-\hat{\lambda}_i(t;\btheta^\star)\right]dt\right]\rightarrow 0\label{2},\\
&\mathbb{E}\left[\frac{1}{T}\int_0^T \left[ \log\left(\frac{\sup_{\btheta\in U_s}{\lambda}_i(t;\btheta)}{\sup_{\btheta\in U_s}\hat{\lambda}_i(t;\btheta)}\right)\right]N_i(dt)\right]\rightarrow 0\label{3},\\
&\mathbb{E}\left[\frac{1}{T}\int_0^T \left[ \log\left(\frac{{\lambda}_i(t;\btheta^\star)}{\hat{\lambda}_i(t;\btheta^\star)}\right)\right]N_i(dt)\right]\rightarrow 0.\label{4}
\end{align}

By the fact that for any $f$ and $g$,
\[
\inf_{x\in U}f(x)-\inf_{x\in U}g(x)=-\sup_{x\in U}(-f(x))+\sup_{x\in U}(-g(x))\leq \sup_{x\in U}(f(x)-g(x)),
\]
 $$\sup_{x\in U}f(x)-\sup_{x\in U}g(x)\leq \sup_{x\in U}(f(x)-g(x)),$$
and the facts that  $\log(1+x)\leq x$ for $x\geq -1$ and  $\lambda_i(t;\btheta)\geq \underline{\mu}$ for any $t$ and $\btheta$, one can see that the following inequality holds:
\[
\begin{split}
 \log\left(\frac{\sup_{\btheta\in U_s}{\lambda}_i(t;\btheta)}{\sup_{\btheta\in U_s}\hat{\lambda}_i(t;\btheta)}\right)&= \log\left(\frac{\sup_{\btheta\in U_s}{\lambda}_i(t;\btheta)-\sup_{\btheta\in U_s}\hat{\lambda}_i(t;\btheta)}{\sup_{\btheta\in U_s}\hat{\lambda}_i(t;\btheta)}+1\right)
 \leq \frac{\sup_{\btheta\in U_s}\left({\lambda}_i(t;\btheta)-\hat{\lambda}_i(t;\btheta)\right)}{\underline{\mu}}.
 \end{split}
\]

Hence (\ref{1}) and (\ref{3}) can be reduced to the following limits as $T\to \infty$,
\begin{align}
&\mathbb{E}\left[\frac{1}{T}\int_0^T \left[ \sup_{\btheta \in U_s}({\lambda}_i(t;\btheta)-\hat{\lambda}_i(t;\btheta))\right]dt\right]\rightarrow 0\label{1'},\\
&\mathbb{E}\left[\frac{1}{T}\int_0^T  \left[ \sup_{\btheta \in U_s}({\lambda}_i(t;\btheta)-\hat{\lambda}_i(t;\btheta))\right]N_i(dt)\right]\rightarrow 0.\label{3'}
\end{align}

Similarly, (\ref{4}) can be reduced to the following limit as $T\rightarrow\infty$:
\begin{equation}
\mathbb{E}\left[\frac{1}{T}\int_0^T  \left[ {\lambda}_i(t;\btheta^\star)-\hat{\lambda}_i(t;\btheta^\star)\right]N_i(dt)\right]\rightarrow 0.\label{4'}
\end{equation}

Therefore,  it suffices to prove that as $T\rightarrow\infty$, for $\btheta\in \Theta$ and in the neighborhood $U\subseteq \Theta$ of $\btheta$, 
\begin{align}
&\mathbb{E}\left[\frac{1}{T}\int_0^T \left[ \sup_{\btheta \in U}({\lambda}_i(t;\btheta)-\hat{\lambda}_i(t;\btheta))\right]dt\right]=\frac{1}{T}\int_0^T \mathbb{E}\left[\sup_{\btheta \in U}({\lambda}_i(t;\btheta)-\hat{\lambda}_i(t;\btheta))\right]dt\rightarrow 0\label{12'},\\
&\mathbb{E}\left[\frac{1}{T}\int_0^T  \left[ \sup_{\btheta \in U}({\lambda}_i(t;\btheta)-\hat{\lambda}_i(t;\btheta))\right]N_i(dt)\right]\rightarrow 0.\label{34'}
\end{align}

To this end, we see by Lemma \ref{stat_moment_inflambda} and by applying Proposition \ref{dN2dt_ogata} to $\xi_i^{(3)}(t)$, 
\[
\begin{split}
\mathbb{E}&\left[\frac{1}{T}\int_0^T  \left[ \sup_{\btheta \in U}({\lambda}_i(t;\btheta)-\hat{\lambda}_i(t;\btheta))\right]N_i(dt)\right]=\mathbb{E}\left[\frac{1}{T}\int_0^T \lambda_i(t;\btheta^\star)\sup_{\btheta \in U}({\lambda}_i(t;\btheta)-\hat{\lambda}_i(t;\btheta))\right]dt\\
&\leq \frac{1}{T}\int_0^T\sqrt{\mathbb{E}[\lambda_i^2(t;\btheta^\star)]}\sqrt{\mathbb{E}[\sup\nolimits_{\btheta \in U}({\lambda}_i(t;\btheta)-\hat{\lambda}_i(t;\btheta))^2]}dt\\
&=\sqrt{\mathbb{E}[\lambda_i^2(0;\btheta^\star)]}\frac{1}{T}\int_0^T\sqrt{\mathbb{E}[\sup\nolimits_{\btheta \in U}({\lambda}_i(t;\btheta)-\hat{\lambda}_i(t;\btheta))^2]}dt.
\end{split}
\]
Moreover, Lemma \ref{newlemma_consist} implies $\mathbb{E}[\lambda_i(0;\btheta^\star)^2]<\infty$, 
and Lemma \ref{approx} shows that $\mathbb{E}[\sup\nolimits_{\btheta \in U}({\lambda}_i(t;\btheta)-\hat{\lambda}_i(t;\btheta))^2]$ is uniformly bounded for all $t\geq 0$ and goes to $0$ as $t\rightarrow\infty$. These results further imply that $\mathbb{E}[\sup\nolimits_{\btheta \in U}|{\lambda}_i(t;\btheta)-\hat{\lambda}_i(t;\btheta)|]$ is uniformly bounded for all $t\geq 0$ and goes to $0$ as $t\rightarrow\infty$. 

The proof is finished by  Proposition \ref{cesaro} and by letting $f(t)=\mathbb{E}[\sup\nolimits_{\btheta \in U}({\lambda}_i(t;\btheta)-\hat{\lambda}_i(t;\btheta))^2]$ for (\ref{34'}) and $f(t)=\mathbb{E}[\sup\nolimits_{\btheta \in U}|{\lambda}_i(t;\btheta)-\hat{\lambda}_i(t;\btheta)|]$ for (\ref{12'}).
\qed

Now,  given a finite amount of time, an optimization problem in general can only be solved with approximation. Therefore it is natural to ask how accurate the computation of $\hat{\btheta}_T$ should be for a given time $T>0$, so that the approximation is also consistent. 
It is easy to see from the above proofs that  consistency still holds if we compute the MLE to $o(T)$ accuracy. 
That is, 
 \begin{corollary}[Consistency of approximate regularized MLEs] \label{parameter_conv_approx}
Under Assumptions \ref{assumption1}, \ref{assumption2}, \ref{assumption4} and \ref{assumption3}, and  assume that $P$ is continuous in $\btheta$. Suppose that $\hat{\btheta}_T^{\text{approx}}$ is obtained so that $\hat{L}^{reg}_T\left(\hat{\btheta}_T^{\text{approx}}\right)>\max\nolimits_{\btheta\in\Theta}\hat{L}^{reg}_T(\btheta)-\epsilon_T$, with $\lim_{T\rightarrow\infty}\epsilon_T/T=0$. Then the approximate regularized MLE\, $\hat{\btheta}^{approx}_T$ converges to ${\btheta}^{\star}$ in probability as $T \rightarrow \infty$.
 \end{corollary}

\section{An alternating minimization algorithm (AA-iPALM)}\label{optimization}
In this section, we propose an Anderson Accelerated inertial Proximal Alternating Linearized Minimization (AA-iPALM)  algorithm   to compute the MLE $\hat{\btheta}_T^{reg}$  in Eqn. (\ref{mle_reg}). We show that AA-iPALM is  guaranteed to  converge globally to stationary points.
The key idea is to view iPALM as a fixed-point iteration. 

We begin by specifying the  settings and the review of   iPALM  and  its convergence behavior (\cite{iPALM}).

\paragraph{Settings.}
In the following, we assume that the triggering function takes a more explicit form of
 \begin{equation}\label{spec_kernel}
 g_{ij}(t;\bigeta)=\sum\nolimits_{m=1}^M\alpha_{ij}^m\phi_m(t;\bbeta_m), \ \ i, j=1, \cdots, K.
 \end{equation}
 Here $\phi_m$ is known as a ``base kernel'' function, $\bigeta=(\balpha,\bbeta)$, $\balpha=(\balpha^1,\cdots,\balpha^M)$ with $\balpha^m=(\alpha_{11}^m,\cdots,\alpha_{KK}^m)$, $\bbeta=(\bbeta_1,\cdots,\bbeta_M)$ with $\bbeta_m\in\mathbb{R}^d$, and $D=K^2+Md$. Note that this specific form has clear intepretability for the parameters and encapsulates most models in the literature. (See \cite{ZZS2013_ode}, ~\cite{DFASS2015} and ~\cite{XFZ2016}.)

We further assume for simplicity that  set $\Theta$ is a Cartesian product of two convex sets $A$ and $B$, with $A\subseteq\mathbb{R}^{K+K^2}$ and $B\subseteq\mathbb{R}^{Md}$ corresponding to $(\bmu,\balpha)$ and $\bbeta$, respectively. 
To be consistent with the non-negativity assumption on $g_{ij}$, we  require  $A\subseteq \mathbb{R}_+^{K+K^2}$, $\phi_m$ taking values in $\mathbb{R}_+$ for all $m=1,\dots,M$, and $\bmu\geq \underline{\mu}$  for all $\bmu\in A_{\bmu}:=\{\bmu
\;|\;\exists\,\balpha, (\bmu,\balpha)\in A\}$. We also assume that $\Theta\subseteq B(0,R)$ for some $R>0$. 

\paragraph{Explicit formula of the log-likelihood function.}
Now denote  $n_j$ as the total number of arrivals of type $j$ during time interval $[0,T]$ for $1 \leq j \leq K$, and $t_s^j$ as  the time of the $s$-th arrival of type $j$ between time $0$ and $T$ for  $s=1,\cdots,n_j$. 
The log-likelihood function in Eqn. (\ref{exact_ll}) can be explicitly written given the specific timestamps of arrivals:
\begin{equation}\label{likelihood_opt}
\begin{split}
\hat{L}_T&(\bmu,\balpha,\bbeta)=-T\sum\nolimits_{i=1}^K\mu_i-\sum\nolimits_{i=1}^K\sum\nolimits_{j=1}^K\sum\nolimits_{m=1}^M\alpha_{i,j}^m\sum\nolimits_{s=t_1^j}^{t_{n_j}^j}\int_s^T\phi_m(t-s;\bbeta_m)dt\\
&+\sum\nolimits_{i=1}^K\sum\nolimits_{t=t_1^i}^{t_{n_i}^i}\log\left(\mu_i+\sum\nolimits_{j=1}^K\sum\nolimits_{m=1}^M\alpha_{ij}^m\sum\nolimits_{\substack{s=t_1^j,\cdots,t_{n_j}^j,~s\leq t}}\phi_m(t-s;\bbeta_m)\right).
\end{split}
\end{equation}
It can be seen from Eqn. (\ref{likelihood_opt}) that $\hat{L}_T$ is convex in $(\bmu,\balpha)$ for any fixed $\bbeta$. This observation, together with the usual low dimensionality of $\bbeta$ in practice, is key for the consideration of alternating minimization type algorithms.


For this section only, we assume additional differentiability conditions on the kernel functions and the regularization term to establish the convergence of the optimization algorithm.
\begin{assumption}[Differentiability]\label{differentiable}
For all $m=1,\dots,M$ and $t\geq 0$,  $\phi_m(t;\bbeta_m)$ is twice continuously differentiable (\textit{i.e.}, $C^2$) in $\bbeta_m$. $P(\btheta)$ is also twice continuously differentiable in $\btheta$.
\end{assumption}

By rewriting the regularized MLE problem in Eqn. (\ref{mle_reg}) as 
\[
\text{minimize}_{(\bmu,\balpha,\bbeta)\in\mathbb{R}^{K+K^2+Md}}~ F(\bmu,\balpha,\bbeta)=\delta_A(\bmu,\balpha)+\delta_B(\bbeta)-\hat{L}_T^{reg}(\bmu,\balpha,\bbeta),
\]
the following Lemma is straightforward: 
\begin{lemma}\label{iPALM-assumptionA}
Under Assumption \ref{differentiable}, $\hat{L}_T^{reg}$ is $C^2$, with $\inf_{\mathbb{R}^{K+K^2+Md}}F(\bmu,\balpha,\bbeta)>-\infty$. In addition, $\delta_A(\bmu,\balpha)$ and $\delta_B(\bbeta)$ are proper and lower semi-continuous functions with $\inf_{\mathbb{R}^{K+K^2}}\delta_A>-\infty$ and $\inf_{\mathbb{R}^{Md}}\delta_B>-\infty$.
\end{lemma}

By Lemma \ref{iPALM-assumptionA} and the compactness of $A$ and $B$, we can define the following three finite constants: \begin{equation}
L_1:=\sup_{(\bmu,\balpha)\in A}~\lambda_{\max}\left(-\nabla^2_{\bmu,\balpha}\hat{L}_T^{reg}(\bmu,\balpha,\bbeta)\right)<\infty,
\end{equation} 
\begin{equation}
L_2:=\sup_{\bbeta\in B}~\lambda_{\max}\left(-\nabla^2_{\bbeta}\hat{L}_T^{reg}(\bmu,\balpha,\bbeta)\right)<\infty,
\end{equation} 
\begin{equation} 
M:=\sup_{(\bmu,\balpha,\bbeta)\in \Theta}~\lambda_{\max}\left(-\nabla^2_{\bmu,\balpha,\bbeta}\hat{L}_T^{reg}(\bmu,\balpha,\bbeta)\right)<\infty.
\end{equation} 

Now, let us review the iPALM algorithm (Algorithm \ref{alg:iPALM}).

 \begin{algorithm}[ht]\label{alg:iPALM}
   \caption{\textbf{inertial Proximal Alternating Linearized Minimization (iPALM)}}
   \label{alg:iPALM}
\begin{algorithmic}[1]
  \STATE \textbf{Input:} Arrival timestamps $t_s^j$ between $[0,T]$, for $s=1,\cdots,n_j$ and $j=1,\cdots,K$. 
  \STATE Construct regularized log-likelihood function $\hat{L}_T^{reg}({\bmu},{\balpha},{\bbeta})$.
  \STATE Specify $\epsilon\in(0,1/2)$, and choose momentum coefficients $\gamma_1^k,~\gamma_2^k\in[0,1-2\epsilon]$. 
  \STATE Set step sizes $\tau_1^k=\dfrac{2(1-\gamma_1^k)}{(1+\gamma_1^k)\bar{L}_1^k}$, 
  $\tau_2^k=\dfrac{2(1-\gamma_2^k)}{(1+\gamma_2^k)\bar{L}_2^k}$, where $\bar{L}_i^k\geq L_i$, $i=1,~2$.
  \STATE Choose initial point ${\btheta}^0=(\bmu^0,\balpha^0,\bbeta^0)\in\Theta$.
   \STATE \textbf{for} {$k=1, 2, \cdots$}
 \[
\left[ \begin{array}{l}
 \bmu^k\\
 \balpha^k
 \end{array}\right]=\Pi_A\left(
 \left[ \begin{array}{l}
 \bmu^{k-1}\\
 \balpha^{k-1}
 \end{array}\right]
 +\tau_1^k\nabla_{\bmu,\balpha}\hat{L}_T^{reg}\left(\bmu^{k-1},\balpha^{k-1},\bbeta^{k-1}\right)
 +\gamma_1^k\left(
 \left[ \begin{array}{l}
 \bmu^{k-1}\\
 \balpha^{k-1}
 \end{array}\right]-
 \left[ \begin{array}{l}
 \bmu^{k-2}\\
 \balpha^{k-2}
 \end{array}\right]
 \right)
 \right),
 \]
 \[
 \bbeta^k=\Pi_B\left(\bbeta^{k-1}+\tau_2^k\nabla_{\bbeta}\hat{L}_T^{reg}(\bmu^k,\balpha^k,\bbeta^{k-1})+\gamma_2^k\left(\bbeta^{k-1}-\bbeta^{k-2}\right)\right).
 \]
 \end{algorithmic}
\end{algorithm}

In Algorithm \ref{alg:iPALM}, since $A$ and $B$ are bounded, the iterations are all bounded.  Lemma \ref{iPALM-assumptionA} then ensures that $\nabla \hat{L}_T$ is globally Lipschitz continuous with uniformly bounded Lipschitz constants
(Remark (iv), Section 2 of \cite{iPALM}). This then validates Assumption A in \cite{iPALM}. Meanwhile, since $\delta_A$ and $\delta_B$ are convex, Remark 4.1 in \cite{iPALM} indicates that a relaxed version of Assumption B (with $1/2$ replaced by $1$)  is also satisfied by the step-sizes and momentum coefficients specified above. We have thus the following guaranteed global convergence result. 

\begin{proposition}\label{glb_conv_set_crit}
Under Assumption \ref{differentiable}, suppose that \,$\btheta^k=(\bmu^k,\balpha^k,\bbeta^k)$ is the iteration sequence generated by iPALM, with $\btheta^0\in \Theta$, then $\btheta^k\in\Theta$ converges to the set of stationary points of $F(\bmu,\balpha,\bbeta)$. More precisely, the limit point set $L:=\{{\btheta}:\;\exists \{k_i\}_{i=1}^{\infty}\text{ s.t. }\lim_{k_i\rightarrow\infty}{\btheta}^{k_i}={\btheta}\}$ is a nonempty, connected and compact subset of the stationary point set $S$ of $F$, where
$S:=\{{\btheta}:\;\nabla_{\bmu,\balpha}\hat{L}_T^{reg}({\bmu,\balpha}) \in N_A({\bmu,\balpha}),~\nabla_{\bbeta}\hat{L}_T^{reg}({\bbeta}) \in N_B({\bbeta})\}$.
Here $N_{C}(\pmb{x})$ is the normal cone of the nonempty convex and closed set $C$ at  point $\pmb{x}$, defined as $N_C(\pmb{x})=\{\pmb{y}\;|\;\pmb{y}^T(\pmb{z}-\pmb{x})\leq 0~\forall \pmb{z}\in C\}$.
\end{proposition}

Hereafter, we fix $\gamma_i^k=\gamma_i$ and $\tau_i^k=\tau_i$, $\bar{L}_i^k=\bar{L}_i$ ($i=1,~2$) to be constant in $k$. Proposition \ref{glb_conv_set_crit} then implies that every fixed point of the iPALM iteration mapping is stationary in the following sense.
\begin{lemma}\label{fix=crit}
Define $H^{iPALM}$ to be the iteration mapping of iPALM, \textit{i.e.}, $H^{iPALM}(\btheta^{k-1},\btheta^{k-2}):=(\btheta^k,\btheta^{k-1})$. Then under Assumption \ref{differentiable}, $H^{iPALM}(\btheta',\btheta)=(\btheta',\btheta)$ if and only if $\btheta=\btheta'$ is a stationary point of  $F$. 
\end{lemma} 

We will see that the above lemma lays the foundation for the incorporation of Anderson acceleration to iPALM in the next section.

\subsection{The AA-iPALM algorithm}

In this section, we propose an AA-iPALM to accelerate the iPALM algorithm. We will show that AA-iPALM inherits the global convergence property of iPALM.

We begin by rephrasing two critical propositions in \cite{iPALM} in terms of the fixed-point  mapping $H^{iPALM}$. 

\begin{proposition}[Propositions 4.3, 4.4 in \cite{iPALM}]\label{ipalm_bounds}
For any $\btheta=(\bmu,\balpha,\bbeta)$ and\, $\btheta'=(\bmu',\balpha',\bbeta')\in \Theta$, denote $H^{iPALM}(\btheta',\btheta)=(\btheta'',\btheta')$ and $\btheta''=(\bmu'',\balpha'',\bbeta'')$. Define $\delta_1=\frac{\gamma_1\bar{L}_1}{2(1-\epsilon-\gamma_1)}$, $\delta_2=\frac{\gamma_2\bar{L}_2}{2(1-\epsilon-\gamma_2)}$, $\rho_1=\frac{\epsilon}{2}\min\{\delta_1,\delta_2\}$, $\rho_2=\sqrt{2}\left(\frac{1-\epsilon}{2\epsilon}\bar{L}_1+\frac{1+\epsilon}{2\epsilon}\bar{L}_2+M\right)+2(\delta_1+\delta_2)$,  where $\epsilon$ is the parameter in Algorithm \ref{alg:iPALM}. Then under Assumption \ref{differentiable}, we have 
\begin{equation}\label{ipalm_bound1}
\begin{split}
\rho_1\|(\btheta''-\btheta',\btheta'-\btheta)\|_2^2&\leq \hat{L}_T^{reg}(\bmu'',\balpha'',\bbeta'')-\hat{L}_T^{reg}(\bmu',\balpha',\bbeta')\\
&+\frac{\delta_1}{2}\left(\|(\bmu'-\bmu,\balpha'-\balpha)\|_2^2-\|(\bmu''-\bmu',\balpha''-\balpha')\|_2^2\right)\\
&+\frac{\delta_2}{2}\left(\|\bbeta'-\bbeta\|_2^2-\|\bbeta''-\bbeta'\|_2^2\right).
\end{split}
\end{equation}
Moreover,  $\bv_{\bmu,\balpha}''\in \partial_{\bmu,\balpha} F(\bmu'',\balpha'',\bbeta'')$, $\bv_{\bbeta}''\in \partial_{\bbeta}F(\bmu'',\balpha'',\bbeta'')$, and 
\begin{equation}\label{ipalm_bound2}
\begin{split}
\|\left(\bv_{\bmu,\balpha}''+\delta_1(\bmu''-\bmu',\balpha''-\balpha'),\bv_{\bbeta}''+\delta_2(\bbeta''-\bbeta')\right)\|_2\leq \rho_2\|(\btheta''-\btheta',\btheta'-\btheta)\|_2.
\end{split}
\end{equation}
Here \[
\bv_{\bmu,\balpha}'':=\nabla_{\bmu,\balpha}\hat{L}_T^{reg}(\bmu',\balpha',\bbeta')-\nabla_{\bmu,\balpha}\hat{L}_T^{reg}(\bmu'',\balpha'',\bbeta'')+\frac{1}{\tau_1}\left((\bmu'-\bmu'',\balpha'-\balpha'')+\gamma_1(\bmu'-\bmu,\balpha'-\balpha)\right),
\]
 \[
\bv_{\bbeta}'':=\nabla_{\bbeta}\hat{L}_T^{reg}(\bmu'',\balpha'',\bbeta')-\nabla_{\bbeta}\hat{L}_T^{reg}(\bmu'',\balpha'',\bbeta'')+\frac{1}{\tau_2}\left(\bbeta'-\bbeta''+\gamma_2(\bbeta'-\bbeta)\right).
\]

\end{proposition}
We are now ready to state the AA-iPALM algorithm (Algorithm \ref{alg:AA-iPALM}). To simplify the notation, we denote $\pmb{u}=(\btheta',\btheta)$, with $\pmb{u}_1=\btheta'$ and $\pmb{u}_2=\btheta$.

\begin{algorithm}[h]\label{alg:AA-iPALM}
   \caption{\textbf{Anderson Accelerated iPALM (AA-iPALM)}}
   \label{alg:AA-iPALM}
\begin{algorithmic}[1]
   \STATE {\bfseries Input:} Arrival time $t_s^j\in[0,T]$ for $s=1,\cdots,n_j$, $j=1,\cdots,K$, step size parameter $\epsilon$, Lipschitz constants (upper bounds) $\bar{L}_1,\bar{L}_2$, momentum coefficients $\gamma_1,\gamma_2$, step sizes $\tau_1,\tau_2$, stabilization constants $\bar{\omega},~\nu\in(0,1)$, safe-guard constants $\delta\geq \max\{\delta_1,\delta_2\}$, $C_1,C_2\geq 1$, max-memory $m>0$. 
   \STATE Choose initial point ${\btheta}^0=(\bmu^0,\balpha^0,\bbeta^0)\in\Theta$, and set $\pmb{u}^0=(\btheta^0,\btheta^0)$.
   \STATE Initialize $H_0=I$, $m_0=0$, 
   and compute $\pmb{u}^1=\tilde{\pmb{u}}^1=H^{iPALM}(\pmb{u}^0)$. Let $\btheta^1=(\bmu^1,\balpha^1,\bbeta^1):=\pmb{u}_1^1$.
   \FOR{$k=1, ~2, ~\dots$}
       \STATE $m_k=m_{k-1}+1$.
       \STATE Compute $\pmb{s}_{k-1}=\tilde{\pmb{u}}^{k}-\pmb{u}^{k-1}$, $\pmb{y}_{k-1}
       =\tilde{\pmb{u}}^{k}-\pmb{u}^{k-1}-\left(H^{iPALM}(\tilde{\pmb{u}}^{k})-H^{iPALM}(\pmb{u}^{k-1})\right)$.
    \STATE Compute $\hat{\pmb{s}}_{k-1}=\pmb{s}_{k-1}-\sum_{j=k-m_k}^{k-2}
    \frac{\hat{\pmb{s}}_j^T\pmb{s}_{k-1}}{\hat{\pmb{s}}_j^T\hat{\pmb{s}}_j}\hat{\pmb{s}}_j$.
    \STATE {\bf If} $m_k=m+1$  or $\|\hat{\pmb{s}}_{k-1}\|_2<\nu\|\pmb{s}_{k-1}\|_2$\,\pmb{:}
      \STATE \hspace{1cm} reset $m_k=0$, $\hat{\pmb{s}}_{k-1}=\pmb{s}_{k-1}$, and $H_{k-1}=I$.
   \STATE Compute $\tilde{\pmb{y}}_{k-1}=\omega_{k-1}\pmb{y}_{k-1}-(1-\omega_{k-1})(\pmb{u}^{k-1}-H^{iPALM}(\pmb{u}^{k-1}))$
   \STATE \hspace{1.7cm} with $\omega_{k-1}=\phi_{\bar{\omega}}\left(\hat{\pmb{s}}_{k-1}^TH_{k-1}\pmb{y}_{k-1}/\|\hat{\pmb{s}}_{k-1}\|^2\right)$.
   \STATE Update $H_k=H_{k-1}+\frac{(\pmb{s}_{k-1}-H_{k-1}\tilde{\pmb{y}}_{k-1})
   \hat{\pmb{s}}_{k-1}^TH_{k-1}}{\hat{\pmb{s}}_{k-1}^TH_{k-1}\tilde{\pmb{y}}_{k-1}}$.
   \STATE Compute $\hat{\pmb{u}}^{k+1}=H^{iPALM}(\pmb{u}^k)$ (iPALM candidate),\\ 
   \hspace{1.7cm}$\tilde{\pmb{u}}^{k+1}=\pmb{u}^k-H_k(\pmb{u}^k-\hat{\pmb{u}}^{k+1})$ (AA candidate).
   \STATE Let $\left(\hat{\bmu}^{k+1},\hat{\balpha}^{k+1},\hat{\bbeta}^{k+1}\right):=\hat{\pmb{u}}^{k+1}_1$, and $\left(\tilde{\bmu}^{k+1},\tilde{\balpha}^{k+1},\tilde{\bbeta}^{k+1}\right):=\tilde{\pmb{u}}^{k+1}_1$. 
   \STATE {\bf If} {\color{black}$\|\nabla\hat{L}_T^{reg}({\bmu}^{k},{\balpha}^{k},{\bbeta}^{k})\|_2\leq C_1\|\hat{\pmb{u}}^{k+1}-\pmb{u}^{k}\|_2$} \& $\left(\tilde{\bmu}^{k+1},\tilde{\balpha}^{k+1},\tilde{\bbeta}^{k+1}\right)\in \Theta$ \\
   \hspace{0.5cm}\& $\|\tilde{\pmb{u}}_2^{k+1}-{\pmb{u}}_2^k\|_2\leq C_2\|\hat{\pmb{u}}_2^{k+1}-{\pmb{u}}_2^k\|_2$ \& $\hat{L}_T^{reg}(\tilde{\bmu}^{k+1},\tilde{\balpha}^{k+1},\tilde{\bbeta}^{k+1})-\hat{L}_T^{reg}(\bmu^k,\balpha^k,\bbeta^k)\geq$\\
   \hspace{1cm}$\frac{\delta+\epsilon\delta}{2}(\|\tilde{\bmu}^{k+1}-\bmu^k\|_2^2+\|\tilde{\balpha}^{k+1}-\balpha^k\|_2^2+\|\tilde{\bbeta}^{k+1}-\bbeta^k\|_2^2)$\,\pmb{:}
set $\pmb{u}^{k+1}=\tilde{\pmb{u}}^{k+1}$. 
   \STATE {\bf else} $\pmb{u}^{k+1}=\hat{\pmb{u}}^{k+1}$.
   \STATE Set $\btheta^{k+1}=(\bmu^{k+1},\balpha^{k+1},\bbeta^{k+1}):=\pmb{u}_1^{k+1}$.
\ENDFOR
\end{algorithmic}
\end{algorithm}
In Algorithm \ref{alg:AA-iPALM}, the function $\phi_{\bar{\omega}}$ in line 11 is defined as 
\begin{equation}\label{phi}
\phi_{\bar{\omega}}(\eta)=
\left\{
\begin{array}{ll}
1 & \text{if $|\eta|\geq \bar{\omega}$}\\
\frac{1-\textbf{sign}(\eta)\bar{\omega}}{1-\eta} 
&\text{if $|\eta|<\bar{\omega}$},
\end{array}
\right.
\end{equation}
where we adopt the convention that $\textbf{sign}(0)=1$. This is exactly the function in \cite{Powell} that ensures non-singularity of the approximate (inverse) Jacobians $H_k$.  
Here, lines 5-9 perform a re-start checking strategy, lines 10-12 perform a Powell-type regularization, both of which introduced in \cite{AA1}. In addition, lines 15-16 execute
a safeguarding strategy with four conditions specially designed for iPALM. 
In particular, the first one ensures that the gradient norm at the current iteration is small, the second one guarantees that the AA candidate for the next iteration lies in the constraint set, the third one makes sure that the AA and iPALM candidates are close to each other, and the last condition is an adaptation of inequality (\ref{ipalm_bound1}). When these four conditions are all satisfied, the AA candidate is chosen and acceleration is achieved. Otherwise, the algorithm proceeds with a vanilla iPALM update. 

Notice that in Algorithm \ref{alg:AA-iPALM}, the iteration variables are $\pmb{u}^k, ~\hat{\pmb{u}}^k,~\tilde{\pmb{u}}^k\in\mathbb{R}^{2(K+D)}$, from which the actual updates $\btheta^k=(\bmu^k,\balpha^k,\bbeta^k)$, the AA candidates $(\hat{\bmu}^k,\hat{\balpha}^k,\hat{\bbeta}^k)$, and the iPALM candidates $(\tilde{\bmu}^k,\tilde{\balpha}^k,\tilde{\bbeta}^k)$ are defined.

The re-start checking strategy and Powell-type regularization ensure that the singular values of the approximate inverse Jacobians are bounded both from above and from below. That is, 
\begin{proposition}[Lemma 3, Corollary 4 in \cite{AA1}]\label{Hkbounds}
In AA-iPALM, we have  $\|H_k\|_2\leq \sigma_H^+:=\left(3\left(\dfrac{1+\bar{\omega}+\nu}{\nu}\right)^m-2\right)^{n-1}/\bar{\omega}^m$ and  $\|H_k^{-1}\|\leq \sigma_H^-:=3\left(\dfrac{1+\bar{\omega}+\nu}{\nu}\right)^m-2$, $\forall k\geq 0$.
\end{proposition} 

Propositions \ref{ipalm_bounds} and \ref{Hkbounds} then guarantee the following global convergence theorem.
\begin{theorem}\label{AA-iPALM-convergence}
Under Assumption \ref{differentiable}, suppose that \,$\btheta^k=(\bmu^k,\balpha^k,\bbeta^k):=\pmb{u}^k_1$ is the iteration sequence generated by AA-iPALM, with $\btheta^0\in \Theta$, then $\btheta^k\in\Theta$ converges to the set of stationary points of $F(\bmu,\balpha,\bbeta)$, in the same sense as in Proposition \ref{glb_conv_set_crit}.
\end{theorem}

A by-product of the proof of Theorem \ref{AA-iPALM-convergence} is the following iterative complexity result.
\begin{theorem}\label{AA-iPALM-complexity}
Under the same condition in Theorem \ref{AA-iPALM-convergence}, we have, for any $K\geq 1$,
\begin{equation*}
\begin{split}
\min_{k\leq K}&~\textbf{dist}(\pmb{0},\partial_{\bmu,\balpha,\bbeta}F(\bmu^k,\balpha^k,\bbeta^k))^2\\
\leq& \dfrac{2\delta^2+2\left(C_2\sigma_H^-\left(\rho_2+\sqrt{2C_1^2+2L_2^2}+\left(M+\delta+\frac{1}{\tau}\right)(1+\sigma_H^+)+\frac{\gamma}{\tau}\right)\right)^2}{\rho_1K}\\
&\times\left(\sup_{(\bmu,\balpha,\bbeta)\in \Theta}\hat{L}_T^{reg}(\bmu,\balpha,\bbeta)-\hat{L}_T^{reg}(\bmu^0,\balpha^0,\bbeta^0)+2{\delta}R^2\right).
\end{split}
\end{equation*}
Here $\tau:=\min\{\tau_1,\tau_2\}$, $\gamma:=\max\{\gamma_1,\gamma_2\}$.
\end{theorem}

Note that AA-iPALM is not restricted to the regularized MLEs for MHPs. It is applicable to general nonconvex optimization problems
as long as the technical assumptions in the above theorems are satisfied.

\section{Experiments}\label{experiments}

\subsection{Experiment Design}
\paragraph{Models.}
Recall that we consider the parameter estimation of  MHPs with a baseline intensity $\pmb{\mu}$ and with the following triggering functions:
 \begin{equation}\label{spec_kernel2}
 g_{ij}(t;\balpha, \bbeta)=\sum\nolimits_{m=1}^M\alpha_{ij}^m\phi_m(t;\bbeta_m), \ \ i, j=1, \cdots, K.
 \end{equation}
 The parameters to be estimated are $\pmb{\theta}=(\pmb{\mu}, \balpha,\bbeta)$, where $\pmb{\mu}=(\mu_1, \cdots, \mu_K)$,  $\balpha=(\balpha^1,\cdots,\balpha^M)$ with $\balpha^m=(\alpha_{11}^m,\cdots,\alpha_{KK}^m)$, 
 and $\bbeta=(\bbeta_1,\cdots,\bbeta_M)$ with $\bbeta_m\in\mathbb{R}^d$.
 
 We first test the algorithms on synthetic datasets, with data generated by \texttt{ticks} (\cite{tick2017}). In particular, we compare the AA-iPALM algorithm with iPALM and PALM in terms of the objective values in the iteration processes. In each experiment, the step sizes $\tau_1,~\tau_2$ and momentum coefficients $\gamma_1,~\gamma_2$ are the same across different algorithms. For AA-iPALM, we choose stabilization constants $\bar{\omega}=\nu=0.1$, safeguard constants $\delta=0.02$, $C_1=C_2=10^8$ (to encourage using AA), and memory size $m=20$. 
 We then apply  AA-iPALM to some real-world datasets, and show that the learned patterns can be well interpreted. 
 
 Throughout this section, we use the quadratic (or Tikhonov) regularization (\cite{YEHK2017}), \textit{i.e.}, $P=C\|\btheta\|_2^2$, where $C$ is some regularization parameter to be specified later.

\subsection{Synthetic experiments}

\subsubsection{Comparison of AA-iPALM, iPALM and PALM}

\begin{figure}
\centering
\includegraphics[width=8.1cm]{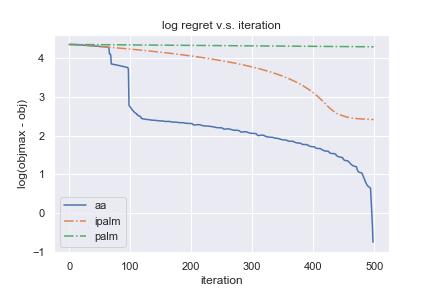}
\includegraphics[width=8.1cm]{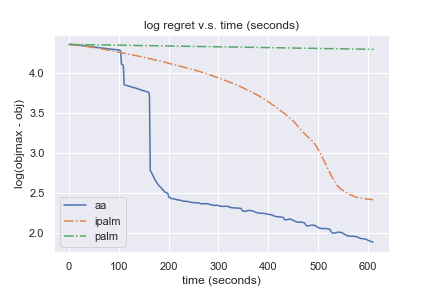}
\caption{Exponential Kernel. Left: log regret v.s. iteration number; Right: log regret v.s. time (seconds).}
\label{exponential_reg_loglikelihood}
\end{figure}

We first consider the model with $K=10$, $M=1$,  and the kernel $\phi_1(t)=e^{-\beta t}$. Here we choose $\bbeta=\beta=0.5$, and generate $\balpha$ by uniformly sampling between $0.001$ and $1$ and then divided by $11$ to ensure stationarity of the MHP. The baseline intensity is similarly generated by uniformly sampling between $0.001$ and $0.1$ and then divided by $2$.  The regularization parameter is set to $C=1$. We choose $\Theta=[lb_{\bmu},ub_{\bmu}]\times[lb_{\balpha},ub_{\balpha}]\times[lb_{\bbeta},ub_{\bbeta}]$, with $lb_{\bmu}=\min_i\mu_i/100$, $ub_{\bmu}=100\max_i\mu_i$, 
$lb_{\balpha}=0$, $ub_{\balpha}=100\max_{i,j}\alpha_{ij}$,
$lb_{\bbeta}=\beta/100$, $ub_{\bbeta}=100\beta$. We initialize $\balpha$ and $\bmu$ with all entries equal to $1$, and $\bbeta=3$. The step sizes are set to $\tau_1=\tau_2=10^{-7}$, and the momentum coefficients are set to $\gamma_1=\gamma_2=0.9$. All three algorithms are run for $500$ iterations. The log regret, \textit{i.e.}, the logarithm of the difference between the maximal regularized log-likelihood value  and the current objective value of the current iteration step, is shown in Figure \ref{exponential_reg_loglikelihood}. 

\begin{figure}
\centering
\includegraphics[width=8.1cm]{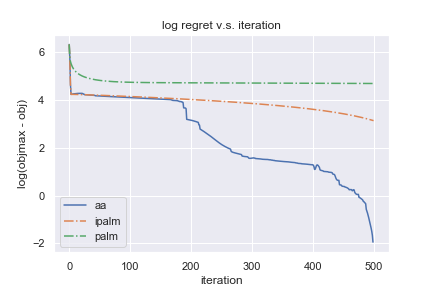}
\includegraphics[width=8.1cm]{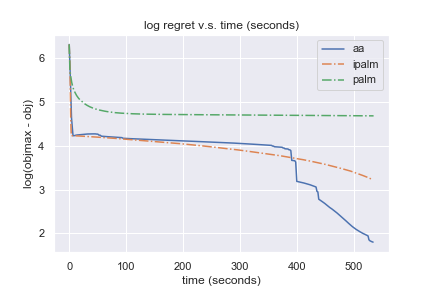}
\caption{Power-law Kernel. Left: log regret v.s. iteration number; Right: log regret v.s. time (seconds).}
\label{pwl_reg_loglikelihood}
\end{figure}

We then consider the power-law kernels, \textit{i.e.}, a model with $K=10$, $M=1$, and $\phi_1(t)=(t+c)^{-\beta t}$, where the cut-off parameter $c$ is chosen to be $0.05$. The problem data is generated in the same way as the exponential kernel example above, except that $\balpha$ is divided by $200$ instead of $11$ to ensure stationarity. The hyper-parameters and initialization are unchanged, except  that $lb_{\bbeta}=\max\{\beta/100,~1.2\}$ as the exponent $\beta$ in power-law kernels is required to be greater than $1$. The corresponding results are shown in Figure \ref{pwl_reg_loglikelihood}. 

Both figures show that  AA-iPALM consistently outperforms iPALM and PALM. Comparing  plots on the left with those on the right, the per iteration cost of AA-iPALM is slightly larger than iPALM, but is compensated by the ``hopping'' acceleration effects. This is especially obvious in the power-law kernel example, in which AA-iPALM mostly follows iPALM for many iterations, and then suddenly ``jumps'' to a much larger value in iteration $200$ or so. In addition, one can see that iPALM consistently outperforms PALM, showing that inertia/momentum does help stabilize the iterations and accelerate the convergence.

\subsection{Applications to real-world data}
We now apply the AA-iPALM algorithm to real-world data. We will use two different sets of data: the Memetracker data and the NASDAQ ITCH data. 
\subsubsection{AA-iPALM for MHPs on the Memetracker dataset}

The Memetracker dataset collects  popular phrases,  the associated articles, their corresponding  url addresses and  publishing time from the internet. It is a popular data source to study the MHPs  (\cite{ZZS2013_ADMM}, \cite{ZZS2013_ode} and \cite{YEHK2017}). Our objective is to study posting activities on social networks, and to analyze the network structure including causality by the MHPs. 

 \paragraph{Data description.}
 We select the most active five news agencies (corresponding to an MHP with $K=5$), and  use the posting data from the ten-day period from April 1st, 2009
to April 10th, 2009.  The top five news agencies are : \texttt{seattle.craigslist.org}, \texttt{chicago.craigslist.org}, \texttt{sfbay.craigslist.org},
 \texttt{blogs.myspace.com}, and  \texttt{sandiego.craigslist.org}. 

\paragraph{Model.} We use the exponential  kernel $
 \phi(t,\bbeta)= \exp(-\bbeta t)$,
 with $M=d=1$, $K=5$, and $\bbeta$ unknown. Here the parameters are $\btheta=(\pmb{\mu}, \balpha,\bbeta)$ with $\pmb{\mu}=(\mu_1, \mu_2,\cdots, \mu_5)$, $\balpha=(\alpha_{11},\cdots,\alpha_{55})$, and $\bbeta=\beta$. This choice of exponential kernel is consistent with the literature. (See \cite{ZZS2013_ADMM}.) 
  
\paragraph{Parameter set-up and initialization.}The regularization parameter is set to $C=0.1$. We choose $\Theta=[lb_{\bmu},ub_{\bmu}]\times[lb_{\balpha},ub_{\balpha}]\times[lb_{\bbeta},ub_{\bbeta}]$, with $lb_{\bmu}=\min_i\mu_i/100$, $ub_{\bmu}=100\max_i\mu_i$, 
$lb_{\balpha}=0$, $ub_{\balpha}=100\max_{i,j}\alpha_{ij}$,
$lb_{\bbeta}=\beta/100$, $ub_{\bbeta}=100\beta$. We initialize $\balpha$ with all entries equal to $0.1$, $\bmu$ with all entries equal to $1$, and $\bbeta=5$. The step sizes are set to $\tau_1=\tau_2=10^{-7}$, and the momentum coefficients are set to $\gamma_1=\gamma_2=0.9$. All three algorithms are run for $500$ iterations.

\begin{figure}[H]
  \begin{subfigure}[b]{0.33\textwidth}
    \includegraphics[width=\textwidth]{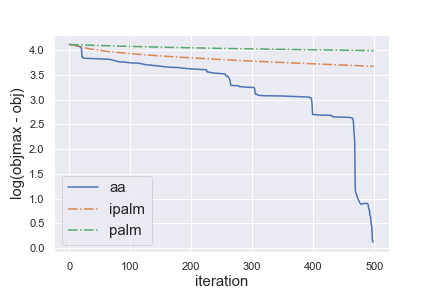}
    \caption{Log regret vs iteration number.}\label{memetrack_iter}
  \end{subfigure}
  \begin{subfigure}[b]{0.33\textwidth}
    \includegraphics[width=\textwidth]{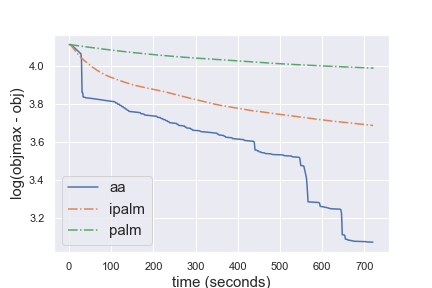}
   \caption{Log regret vs time (seconds).}\label{memetrack_time}
  \end{subfigure}
  \begin{subfigure}[b]{0.33\textwidth}
    \includegraphics[width=\textwidth]{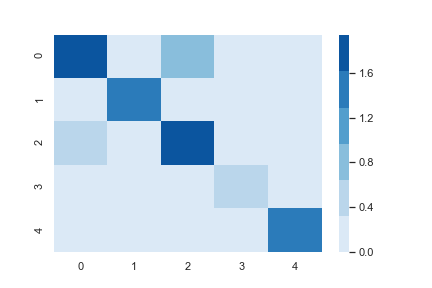}
    \caption{$\pmb{\hat{\alpha}}$ matrix.} \label{memetrack_alpha}
  \end{subfigure}
 \caption{Performance on Memetracker with exponential distribution.}\label{figure: memetracker}

\end{figure}

 \paragraph{Results.} Looking at the convergence curves, one can see that AA-iPALM again consistently outperforms iPALM and PALM. The estimated $\balpha$-matrix is sparse and full rank. The concentration is on the main diagonal (Figure \ref{memetrack_alpha}). It  suggests that people from one area show little interest in rental postings from other areas, with the only exception being the upper middle portion in Figure \ref{memetrack_alpha}, which shows the mutual interest between the San Francisco bay area and the Seattle area. This may indicate active reallocations
of  engineers from tech companies in  these two regions. {The estimated $\beta$ value $\hat{\beta}=2.58$ suggests a moderate decay in postings. } 



\subsubsection{AA-iPALM for MHPs on NASDAQ ITCH dataset}

In financial market, a limit order book (LOB) is used to record and match buying/bid and ask/selling orders of three types: limit, market, and cancellation. 
 A limit order is an order to trade a certain amount
of securities (stocks, futures, etc.) at a specified price. 
A market order is an order to buy/sell a certain amount of the equity at the best available price in the LOB;
it is then matched with the best available price and a trade occurs immediately and the
LOB is updated accordingly. A limit order stays in the LOB until it is executed against a
market order or until it is canceled, and cancellation is allowed at any time. (For more background  on LOB, see  \cite{LL2013}, \cite{CJP2015}, and \cite{GSLW2017}).

\paragraph{Data. }

The NASDAQ ITCH message data from {\it LobsterData} contains all the updates of the LOB from different types of events. The types of events include cancellations, limit orders and market orders. For each update, it contains the information of ``time stamp'', ``event type'', ``Order ID'', ``size'', ``price'', and ``direction''.
 Our focus is to use MHPs to analyze two aspects in LOB and the high frequency trading:  cancellation activities and the herd behavior.

\paragraph{Data description.}
We choose all the order information for  Google (GOOG) and Apple (APPL) 
 from 1:30pm to 3:30pm EST in June 21st, 2012 June 21. The orders include six different types:
 limit orders on the bid/buy side ($L^b$), limit orders on the ask/sell side ($L^a$),  market orders on the bid/buy side ($M^b$), market orders on the ask/sell side ($M^a$), cancellations on the bid/buy side ($C^b$), and cancellations on the ask/sell side ($C^a$). These correspond to an MHP with $K=6$.

\paragraph{Model.} In this experiment, we use the power-law kernel $\phi(t)=(t+c)^{-\beta}$
with $M=d=1$, $K=6$,  a cut-off $c= 0.05$, and $\bbeta$ unknown. Here the parameters are $\btheta=(\pmb{\mu}, \balpha,\bbeta)$ with $\pmb{\mu}=(\mu_1, \mu_2,\cdots, \mu_{6})$, $\balpha=(\alpha_{11},\cdots,\alpha_{66})$, and $\bbeta=\beta$. The choice of power-law kernel is consistent with the literature. (See \cite{EBK2012} and \cite{BJM2016}.)

\begin{figure}[ht]
\begin{subfigure}[h]{0.5\linewidth}
\centering
\includegraphics[width=0.9\linewidth]{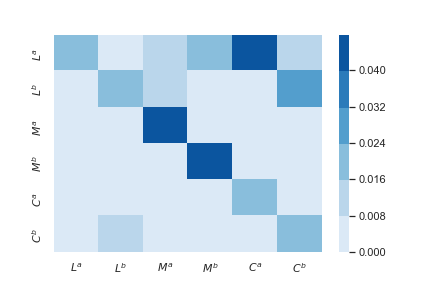}
\caption{Google.}
\end{subfigure}
\hfill
\begin{subfigure}[h]{0.5\linewidth}
\centering
\includegraphics[width=0.9\linewidth]{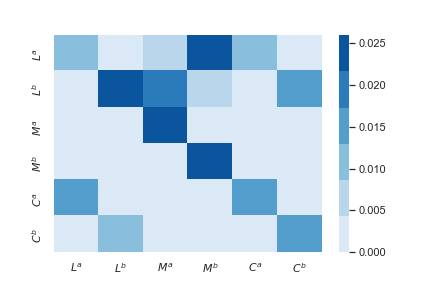}
\caption{Apple.}
\end{subfigure}%
\caption{Performance on NASDAQ data with power-law distribution.}\label{figure: alphamatrix}
\end{figure}

\paragraph{Parameter set-up and initialization.}The regularization parameter is set to $C=0.1$. We choose $\Theta=[lb_{\bmu},ub_{\bmu}]\times[lb_{\balpha},ub_{\balpha}]\times[lb_{\bbeta},ub_{\bbeta}]$, with $lb_{\bmu}=\min_i\mu_i/100$, $ub_{\bmu}=100\max_i\mu_i$, 
$lb_{\balpha}=0$, $ub_{\balpha}=100\max_{i,j}\alpha_{ij}$,
$lb_{\bbeta}=\beta/100$, $ub_{\bbeta}=100\beta$. We initialize $\balpha$ and $\bmu$ with all entries equal to $1$ and $\bbeta=5$. The step sizes are set to $\tau_1=\tau_2=10^{-7}$, and the momentum coefficients are set to $\gamma_1=\gamma_2=0.9$. The cut-off parameter is set as 0.005. The AA-iPALM is run for $500$ iterations.

\paragraph{Results.}
The results for Google and Apple are similar in terms of both the $\pmb{\alpha}$-matrix and value of $\pmb{\mu}$. (Figure \ref{figure: alphamatrix}  and Table~\ref{table:mu_learnt}.) This is no surprising as both companies are leaders in the technology sector and are the most two valuable companies by the market capitalization.
The estimated  values of $\beta$: $\hat{\beta}=1.36$ for Google and $\hat{\beta}=1.44$ for Apple, are also consistent  with those reported by \cite{EBK2012} and \cite{BJM2016}.
Moreover, 

\begin{itemize}
\item non-zero  positive entries on the main diagonal in the $\pmb{\alpha}$-matrix show the self-exciting trading behavior. This is consistent with the well-recognized persistence of order flows, a result of the standard order execution practice to  split  meta-orders into a sequence of smaller orders;
\item non-zero positive entries for $\alpha_{L^a,C^a}$, $\alpha_{L^b,C^b}$, $\alpha_{C^a,L^a}$, and $\alpha_{C^b,L^b}$ are consistent with the well-known high cancellation rate in the high frequency trading.
Most (buy/bid and ask/sell) orders are canceled within seconds of their submissions
 (\cite{HH2007}, \cite{EBK2012}, \cite{BC2013}, \cite{BSKB2015}, and \cite{BJM2016});
\item non-zero entries for $\alpha_{L^a,M^b}$ and $\alpha_{L^b,M^a}$
highlight the herd behavior, termed as the ``same direction (buy or sell) of trading'' in LOB. (\cite{EBK2012} and \cite{BJM2016}.)
 \end{itemize}

 \begin{table}[H]
\centering \small{
\begin{tabular}{ | l | r  r r r r r|}
    \hline
   &$\mu_{1}$&$\mu_{2}$&$\mu_{3}$&$\mu_{4}$&$\mu_{5}$&$\mu_{6}$ \\ \hline \hline
Google $\pmb{\hat{\mu}}$& 0.32& 0.37& 0.23&0.22 &0.32 &0.33\\
Apple $\pmb{\hat{\mu}}$&0.28& 0.25& 0.17& 0.19&0.24& 0.21\\
\hline
   \end{tabular}
    \caption{Estimated baseline intensities $\pmb{\hat{\mu}}$.}
    \label{table:mu_learnt}}
\end{table}



\bibliographystyle{apa}
\bibliography{ref_or}




\newpage
\ECSwitch


\ECHead{Appendices}

In this Appendix, to simplify notations, we extend $g_{ij}$ to be defined on $\mathbb{R}\times\Theta_{\bigeta}$ by setting $g_{ij}(t;\bigeta)=0$ for any $t<0$ and $\bigeta\in\Theta_{\bigeta}$.

\section{Cumulant density formula via family/category trees} Given any time vector $\pmb{t}=(t_1,\dots,t_n)$ and type vector $\pmb{i}=(i_1,\dots,i_n)\in\{1,\dots,K\}^n$,  the $n$-th order cumulant density of the MHP is define as $k^{\pmb{i}}(\pmb{t}):=\dfrac{k(N^{i_1}(dt_1),\dots,N^{i_n}(dt_n))}{dt_1\dots dt_n}$, where $k(X_1,\dots,X_n)$ is the $n$-th order cumulant of random variables $X_1,\dots,X_k$. For example, $k(X)=\mathbb{E}[X]$ and $k(X_1,X_2)=Cov(X_1,X_2)$. Since moments can be expressed as summations and products of cumulants, the desired moments bounds can be derived from bounds on the cumulants.

Under Assumptions \ref{assumption1} and \ref{assumption2}, it follows from \cite{JHR2015} that $$k^{\pmb{i}}(\pmb{t})d\pmb{t}=\mathbb{P}\left(E_{\pmb{t}}^{\pmb{i}}\cap C_{\pmb{t}}^{\pmb{i}}\right),$$ where 
\[
\begin{split}
&E_{\pmb{t}}^{\pmb{i}}=\{\forall k=1,\dots,n, ~\mbox{there is a type $i_k$ event at time $t_k$}\},\\
&C_{\pmb{t}}^{\pmb{i}}=\{\mbox{$\exists$ cluster $C$ such that, $\forall k=1,\dots,n$, $t_k \in C$}\},
\end{split}
\]
and $E_{\pmb{t}}^{\pmb{i}}\cap C_{\pmb{t}}^{\pmb{i}}$ means that there is a type $i_k$ event at time $t_k$ and that all of these events are descendants from a common immigrant. To compute this probability,  notice that all possible branching trees in $E_{\pmb{t}}^{\pmb{i}}\cap C_{\pmb{t}}^{\pmb{i}}$ can be grouped into a finite number of categories. Therefore $\mathbb{P}\left(E_{\pmb{t}}^{\pmb{i}}\cap C_{\pmb{t}}^{\pmb{i}}\right)$ reduces to the sum of the probabilities of these categories. 

One can define the notion of \textit{nearest common ancestor}, where $u$ is called the nearest common ancestor of $v_1,\dots,v_k$ if each $v_i$ ($i=1,\dots,k$) is either equal to $u$ or is a descendant of $u$, and $u$ is the node with the largest time stamp that has this property. For each rooted branching tree in  $E_{\pmb{t}}^{\pmb{i}}\cap C_{\pmb{t}}^{\pmb{i}}$ with root $x$ (immigrant), we keep the root $x$, the type $i_k$ event at time $t_k$ for $k=1,\dots,K$ (represented as $t_k$ for short), and the nearest common ancestors of all the subsets of $\{t_1,\dots,t_n\}$ (the set of which denoted as $A$), and contract the edges which have at least one end point not in $A$. After this operation, each edge in the resulting family/category tree can represent arbitrary number of generations. The idea of the above operation is marginalization, \textit{i.e.}, integration of the joint probabilities over the intermediate generations over time and types.

\begin{figure}[ht]
\centering
\includegraphics[width=14cm]{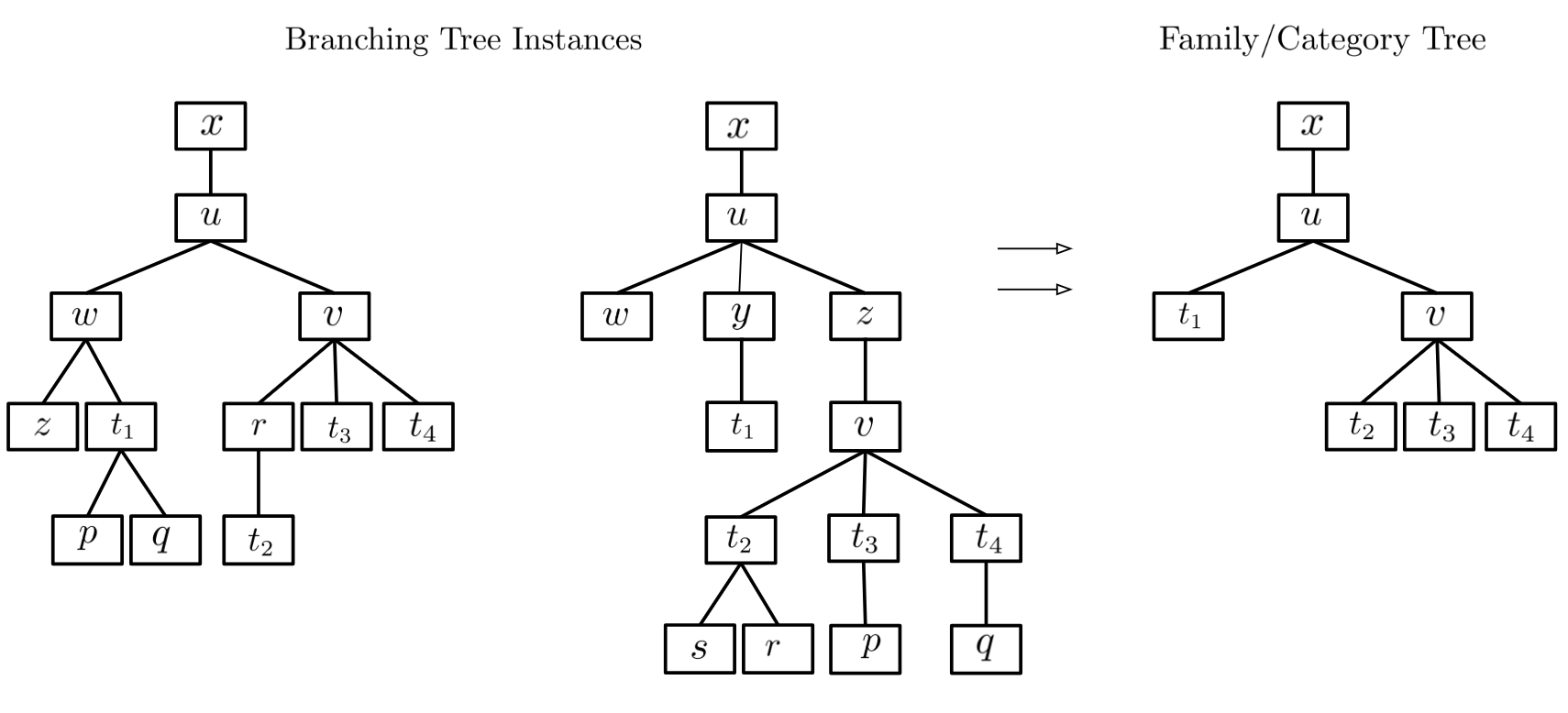}
\caption{Left: branching trees; Right: family/category trees}
\label{family-tree-ex}
\end{figure}

One example is illustrated in Figure \ref{family-tree-ex}. The left-hand side lists two realizations of branching trees, both of which reduce to the same family/category tree on the right-hand side. Notice that nodes $w$, $y$ and $z$ in the branching trees are removed in the family/category tree because none of them is in the set $A=\{u, v\}$ (the nearest common ancestor set).

To compute the probability of each family/category tree, we first compute the product of the conditional probability densities along the edges given the types and event time stamps of all the nodes. We then sum over all possible types and integrate over all possible time stamps.

More precisely, define $R_t^{ij}:=\left[\sum_{n\geq 0}\pmb{G}_{ij}^{\star n}(t)\right]_{ij}$, where $\pmb{G}(t)=[g_{ij}(t;\pmb{\eta}^\star)]_{n\times n}$ is the matrix of triggering functions, and $\pmb{G}^{\star n}(t)$ is the $n$-th (self) convolution of $\pmb{G}$ defined recursively as $\pmb{G}^{\star 0}(t)=\pmb{I}\delta(t)$, $\pmb{G}^{\star n}(t)=\int_{-\infty}^t\pmb{G}^{\star (n-1)}(t-s)\pmb{G}(s)ds$. Here $\delta(t)$ is the Dirac $\delta$ function, and $\int_{-\infty}^{\infty}\pmb{G}(t)dt=\pmb{G}$, where $\pmb{G}$ is the matrix in Assumption \ref{assumption2}. It is shown in \cite{JHR2015} that
\[
R_t^{ij}dt=\mathbb{P}(\mbox{type $j$ event at $0$ causes type $i$ event at $t$}).
\]

Suppose that $u$ is of type $m$, then the probability density along edge $(x,u)$ with the time stamp and type of $x$ marginalized is $\bar{\lambda}_m$, and the probability density along the edge $(u,t_k)$ is $R_{t_k-u}^{i_km}$. Finally, since two nearest common ancestors can not be identical, the probability density along edge $(u,v)$ which connects two (different) nearest common ancestors with types $j$ and $i$ is $\Psi_{v-u}^{ij}:=R_{v-u}^{ij}-\delta_{ij}\delta(v-u)=\left[\sum_{n\geq 1}\pmb{G}^{\star n}(v-u)\right]_{ij}$.

For example, the probability corresponding to the family/category tree in Figure \ref{family-tree-ex} is equal to $\sum_{j_1,j_2=1}^K\bar{\lambda}_{j_2}\int_{\mathbb{R}}R_{t_1-u}^{i_1j_2}\left(\int_{\mathbb{R}}R_{t_2-v}^{i_2j_1}R_{t_3-v}^{i_3j_1}R_{t_4-v}^{i_4j_1}\Psi_{v-u}^{j_1j_2}dv\right)du$.

The above discussion results in the cumulant computation algorithm at the end of Section III in \cite{JHR2015}.

\section{Additional proofs for Section \ref{asymptotic} }

\subsection{Proof of Lemma \ref{moment}.}

\begin{figure}[ht]
\begin{center}
\includegraphics[width=13cm]{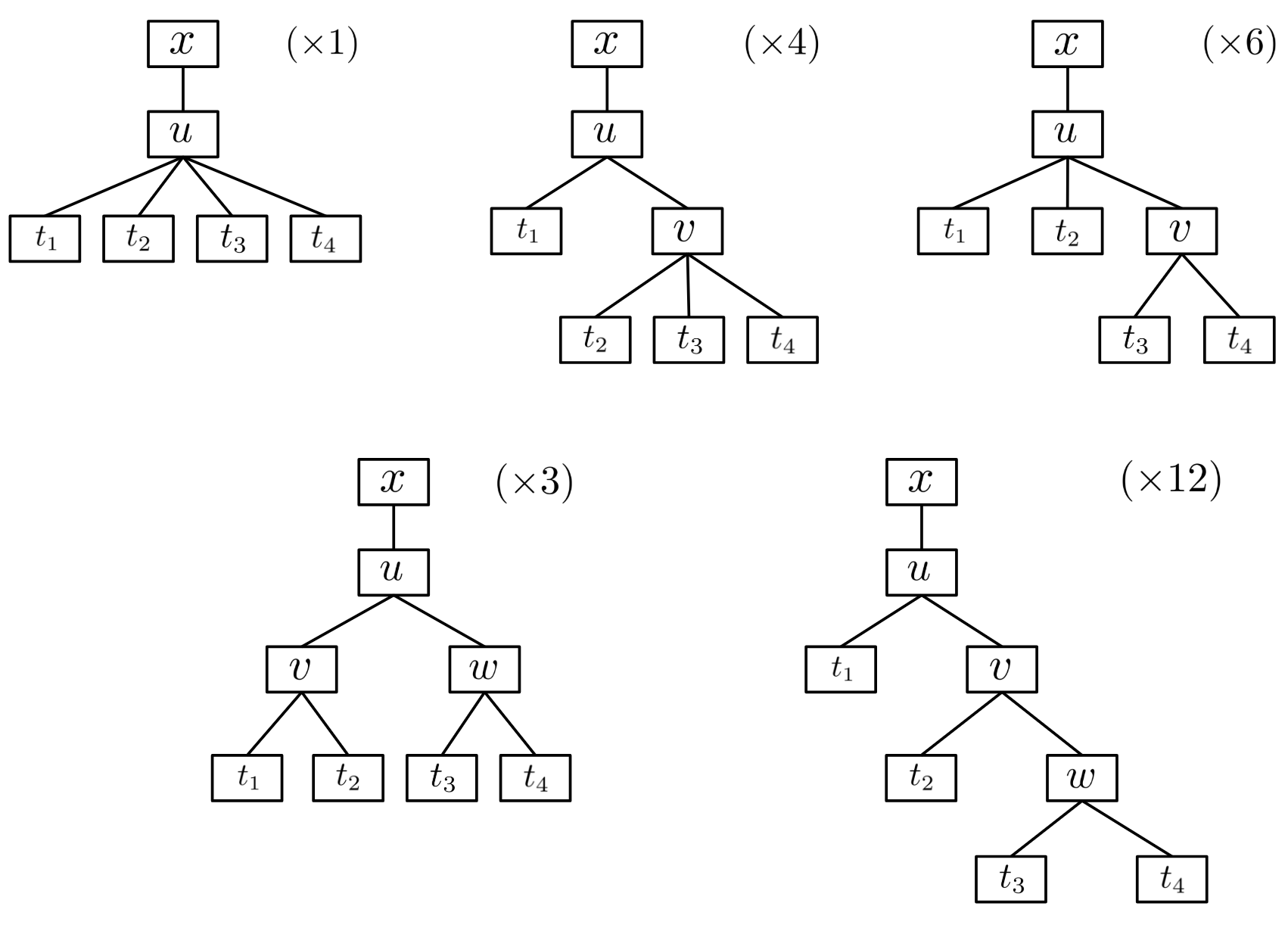}
\caption{A schematic representation of all possible family category trees. The $\times n$ notation at the top right of each tree indicates that by combination and permutation of $(t_1,t_2,t_3,t_4)$, there are $n$ different versions of trees with the same type.}\label{tree}
\end{center}
\end{figure}

The idea is to connect MHPs with the Poisson clustering process outlined above.
For the $4$-th order cumulant, there are 26 family/category trees, which can be further grouped into $5$ generic types by symmetry, as listed in Figure \ref{tree}. 


Hence for $i=j=k=l=i_0$,
\[
\begin{split}
k_4(\pmb{t})&:=k^{i_0i_0i_0i_0}(t_1,t_2,t_3,t_4)=\sum_{j=1}^K\bar{\lambda}_j\int_{\mathbb{R}}\prod_{i=1}^4R_{t_i-u}^{i_0j}du\\
&+\sum_{i=1}^4\sum_{j_1,j_2=1}^K\bar{\lambda}_{j_2}\int_{\mathbb{R}}R_{t_i-u}^{i_0j_2}\left(\int_{\mathbb{R}}\prod_{j\neq i}R_{t_j-v}^{i_0j_1}\Psi_{v-u}^{j_1j_2}dv\right)du\\
&+\sum_{1\leq i_1<i_2\leq 4}\sum_{j_1,j_2=1}^K\bar{\lambda}_{j_2}\int_{\mathbb{R}}R_{t_{i_1}-u}^{i_0j_2}R_{t_{i_2}-u}^{i_0j_2}\left(\int_{\mathbb{R}}\prod_{j\neq i_1,i_2}R_{t_j-v}^{i_0j_1}\Psi_{v-u}^{j_1j_2}dv\right)du\\
&+\sum_{i=2,3,4}\sum_{j_1,j_2,j_3=1}^K\bar{\lambda}_{j_3}\int_{\mathbb{R}}\left(\int_{\mathbb{R}}R_{t_1-v}^{i_0j_2}R_{t_i-v}^{i_0j_2}\Psi_{v-u}^{j_2j_3}dv\right)\left(\int_{\mathbb{R}}\prod_{j\neq 1,i}R_{t_j-w}^{i_0j_1}\Psi_{w-u}^{j_1j_3}dw \right) du\\
&+\sum_{i=1}^4\sum_{j_1,j_2,j_3=1}^K\bar{\lambda}_{j_3}\int_{\mathbb{R}}R_{t_i-u}^{i_0j_3}\sum_{j\neq i}\left(\int_{\mathbb{R}}R_{t_j-v}^{i_0j_2}\Psi_{v-u}^{j_2j_3}\left(\int_{\mathbb{R}}\prod_{l\neq i,j}R_{t_l-w}^{i_0j_1}\Psi_{w-v}^{j_1j_2}dw\right)dv\right)du,
\end{split}
\]
where $\pmb{t}:=(t_1,t_2,t_3,t_4)$, and $k^{ijkl}(t_1,t_2,t_3,t_4)$ is the $4$-th order cumulant of $\pmb{N}(dt)$, as defined in \cite{JHR2015}.

Noticing that $R_t^{ij}$ and $\Psi_t^{ij}$ are both non-negative by definition, $\int_{\mathbb{R}}R_t^{ij}dt=[(\pmb{I}-\pmb{G})^{-1}]_{ij}:=A_{ij}$, and $\int_{\mathbb{R}}\Psi_t^{ij}dt=[\pmb{G}(\pmb{I}-\pmb{G})^{-1}]_{ij}:=B_{ij}$, by integrating from $0$ to $h$ for an arbitrary $h>0$, 
\[
\begin{split}
K_4:=&\int_0^h\int_0^h\int_0^h\int_0^hk_4(t_t,t_2,t_3,t_4)dt_1dt_2dt_3dt_4\\
\leq& \sum_{j=1}^K\bar{\lambda}_jA_{i_0j}^4h+4\sum_{j_1,j_2=1}^K\bar{\lambda}_{j_2}A_{i_0j_2}A_{i_0j_1}^3B_{j_1j_2}h^2\\
&+6\sum_{j_1,j_2=1}^K\bar{\lambda}_{j_2}A_{i_0j_2}^2A_{i_0j_1}^2B_{j_1j_2}h^2+3\sum_{j_1,j_2,j_3=1}^K\bar{\lambda}_{j_3}A_{i_0j_2}^2B_{j_2j_3}A_{i_0j_1}^2B_{j_1j_3}h^3\\
&+12\sum_{j_1,j_2,j_3=1}^K\bar{\lambda}_{j_3}A_{i_0j_3}A_{i_0j_2}B_{j_2j_3}A_{i_0j_1}^2B_{j_1j_2}h^3=O(h^3).
\end{split}
\]
Notice that here we leave out one of $t_1,\dots,t_4$ to maintain dependence on $h$. In this way, we have
\[
\begin{split}
&\int_0^h\int_0^h\int_0^h\int_0^h\int_{\mathbb{R}}R_{t_1-u}^{i_0j}R_{t_2-u}^{i_0j}R_{t_3-u}^{i_0j}R_{t_4-u}^{i_0j}dudt_1dt_2dt_3dt_4\\
&\leq \int_0^h\left(\int_{\mathbb{R}}\left(\int_{\mathbb{R}}R_{t_1-u}^{i_0j}dt_1\right)\left(\int_{\mathbb{R}}R_{t_2-u}^{i_0j}dt_2\right)\left(\int_{\mathbb{R}}R_{t_3-u}^{i_0j}dt_3\right)R_{t_4-u}^{i_0j}du\right)dt_4\\
&=A_{i_0j}^3\int_0^h\left(\int_{\mathbb{R}}R_{t_4-u}^{i_0j}du\right)dt_4=A_{i_0j}^4h.
\end{split}
\]
The rest are similar.

Similarly, by Eqn. (14), (37) and (39) in \cite{JHR2015}, and by setting all the indices to $i_0$, we obtain
$K_1:=\int_0^hk^{i_0}(t)dt=\bar{\lambda}_{i_0}h=O(h)$, $K_2:=\int_0^h\int_0^hk^{i_0i_0}(t_1,t_2)dt_1dt_2=O(h)$,
$K_3:=\int_0^h\int_0^h\int_0^hk^{i_0i_0i_0}(t_1,t_2,t_3)dt_1dt_2dt_3=O(h^2)$.

Moreover, using the relation between cumulants and moments, we see that for any $i_0=1,\dots,K$,
$\mathbb{E}[(N_{i_0}(0,h])^4]=K_4+4K_3K_1+3K_2^2+6K_2K_1^2+K_1^4
=\bar{\lambda}_{i_0}^4h^4+O(h^3)$, $\mathbb{E}[N_{i_0}(0,h]]=K_1=\bar{\lambda}_{i_0}h$, $\mathbb{E}[(N_{i_0}(0,h])^2]=K_2+K_1^2=\bar{\lambda}_{i_0}^2h^2+O(h)$, $\mathbb{E}[(N_{i_0}(0,h])^3]=K_3+3K_2K_1+K_1^3=\bar{\lambda}_{i_0}^3h^3+O(h^2)$. Thus, for any $i_0=1,\dots,K$, $\mathbb{E}[(N_{i_0}(0,h]-\bar{\lambda}_ih)^4]=\mathbb{E}[|N_{i_0}(0,h]|^4]-4\mathbb{E}[(N_{i_0}(0,h])^3]\bar{\lambda}_{i_0}h
+6\mathbb{E}[(N_{i_0}(0,h])^2]\bar{\lambda}_{i_0}^2h^2-4\mathbb{E}[(N_{i_0}(0,h])]\bar{\lambda}_{i_0}^3h^3+\bar{\lambda}_{i_0}^4h^4
=O(h^3).
$

Finally, the proof is complete by taking a maximum over $i_0=1,\dots,K$.
\qed

\subsection{Proof of Lemma \ref{identifiability}.} 
The proof again relies on the explicit formula for covariance density in \cite{JHR2015}.
The key is to notice that if $\lambda_i(0;\btheta)=\lambda_i(0;\btheta')$ a.s., then the expectation and variance of $\lambda_i(0;\btheta)-\lambda_i(0;\btheta')$ should both be equal to $0$. Expanding using Eqn. (\ref{trueint}), and defining $h_{ij}(t):=g_{ij}(t;\pmb{\eta})-g_{ij}(t;\pmb{\eta}')$, we have
\begin{equation}\label{first-order}
0=\mathbb{E}\left[\lambda_i(0;\btheta)-\lambda_i(0;\btheta')\right]=\mu_i-\mu_i'+\sum\nolimits_{j=1}^K\mathbb{E}\left[\int_{-\infty}^0h_{ij}(-s)N_j(ds)\right],
\end{equation}
and
\[
\begin{split}
0=\text{Var}&\left(\lambda_i(0;\btheta)-\lambda_i(0;\btheta')\right)=\mathbb{E}\left[\left(\lambda_i(0;\btheta)-\lambda_i(0;\btheta')\right)^2\right]\\
=&\mathbb{E}\left[\left(\mu_i-\mu_i'+\sum\nolimits_{j=1}^K\int_{-\infty}^0h_{ij}(-s)N_j(ds)\right)^2\right]\\
=&(\mu_i-\mu_i')^2+2(\mu_i-\mu_i')\sum\nolimits_{j=1}^K\mathbb{E}\left[\int_{-\infty}^0h_{ij}(-s)N_j(ds)\right]\\
&+\sum\nolimits_{j_1=1}^K\sum\nolimits_{j_2=1}^K\int_{-\infty}^0\int_{-\infty}^0h_{ij_1}(-s_1)h_{ij_2}(-s_2)\mathbb{E}[N_{j_1}(ds_1)N_{j_2}(ds_2)]\\
=&\sum\nolimits_{j_1=1}^K\sum\nolimits_{j_2=1}^K\int_{-\infty}^0\int_{-\infty}^0h_{ij_1}(-s_1)h_{ij_2}(-s_2)\mathbb{E}[N_{j_1}(ds_1)N_{j_2}(ds_2)]-(\mu_i-\mu_i')^2\\
=&\sum\nolimits_{j_1=1}^K\sum\nolimits_{j_2=1}^K\int_{-\infty}^0\int_{-\infty}^0h_{ij_1}(-s_1)h_{ij_2}(-s_2)\left(\mathbb{E}[N_{j_1}(ds_1)N_{j_2}(ds_2)]-\mathbb{E}[N_{j_1}(ds_1)]\mathbb{E}[N_{j_2}(ds_2)]\right)\\
=&\sum\nolimits_{j_1=1}^K\sum\nolimits_{j_2=1}^K\int_{-\infty}^0\int_{-\infty}^0h_{ij_1}(-s_1)h_{ij_2}(-s_2)k^{j_1j_2}(s_1,s_2)ds_1ds_2.
\end{split}
\]

By Eqn. (14) and (37) in \cite{JHR2015}, this means that for any $i=1,\dots,K$,
\[
\begin{split}
0&=\sum\nolimits_{j_1=1}^K\sum\nolimits_{j_2=1}^K\int_{-\infty}^0\int_{-\infty}^0h_{ij_1}(-s_1)h_{ij_2}(-s_2)\sum\nolimits_{m=1}^K\bar{\lambda}_m\int_{\mathbb{R}}R_{s_1-u}^{j_1m}R_{s_2-u}^{j_2m}duds_1ds_2\\
&=\sum\nolimits_{m=1}^K\bar{\lambda}_m\int_{\mathbb{R}}\left(\sum\nolimits_{j_1=1}^K\sum\nolimits_{j_2=1}^K\int_{-\infty}^0\int_{-\infty}^0h_{ij_1}(-s_1)h_{ij_2}(-s_2)R_{s_1-u}^{j_1m}R_{s_2-u}^{j_2m}ds_1ds_2\right)du\\
&=\sum\nolimits_{m=1}^K\bar{\lambda}_m\int_{\mathbb{R}}\left(\sum\nolimits_{j=1}^K\int_{-\infty}^0h_{ij}(-s)R_{s-u}^{jm}ds\right)^2du\\
&=\sum\nolimits_{m=1}^K\bar{\lambda}_m\int_{\mathbb{R}}\left(\sum\nolimits_{j=1}^K\int_0^\infty h_{ij}(s)R_{-u-s}^{jm}ds\right)^2du\\
&=\sum\nolimits_{m=1}^K\bar{\lambda}_m\int_{\mathbb{R}}\left(\sum\nolimits_{j=1}^K(h_{ij}\star R^{jm})_{-u}\right)^2du.
\end{split}
\]
Since $\bar{\lambda}_i>0$ for all $i=1,\dots,K$, this implies that $\sum\nolimits_{j=1}^K(h_{ij}\star R^{jm})_{-u}\equiv 0$ a.e. for all $u\in\mathbb{R}$ and $i,m=1,\dots,K$. Taking the Laplace transform evaluated at $t$, we see that
\[
\sum\nolimits_{j=1}^K\hat{h}_{ij}(t)A_{jm}(t)=0 ~ \text{a.e.}, \quad \mbox{for all}~ i,m=1,\dots,K,
\]
where $A_{ij}(t):=[(I-\hat{\pmb{G}}(t))^{-1}]_{ij}=\hat{R}^{ij}(t)$, in which $\pmb{G}(t):=[g_{ij}(t;\pmb{\eta}^\star)]_{K\times K}$. 

Rewriting in a matrix form, we see that $\hat{\pmb{H}}(t)(I-\hat{\pmb{G}}(t))^{-1}=0$ a.e., where $\hat{\pmb{H}}(t)=[\hat{h}_{ij}(t)]_{K\times k}$. This implies that $\hat{\pmb{H}}(t)\equiv 0$ a.e. for all $t\geq 0$, which holds iff $\pmb{H}(t):=[h_{ij}(t)]_{K\times K}\equiv 0$ a.e. Hence  $g_{ij}(t;\pmb{\eta})=g_{ij}(t;\pmb{\eta}')$ a.e. {in $t$}, and $\pmb{\eta}=\pmb{\eta}'$ by {Assumption \ref{assumption4}} in Section \ref{prob_set}. Finally, plugging in Eqn. (\ref{first-order}), we see that $\mu_i=\mu_i'$ for $i=1,\dots,K$. 
\qed

\subsection{Proof of Lemma \ref{newlemma_consist} and Lemma \ref{approx}.} The proof is based on Lemma \ref{uniform_increment} regarding the uniform moment bounds of successive increment of $\pmb{N}$.

Suppose that $\{t_k\}_{k\geq 0}$ is the sequence in Assumption \ref{assumption3}. For notational simplicity, we use the shorthand
$C_i:=C_i^{(0)}$. Then since $\log(x)\leq x-1\leq x$ for all $x>0$ and $\Theta\in B(0,R)$, 
\[
\begin{split}
&\sup_{\pmb{\theta}'\in \Theta}\lambda_i(0;\pmb{\theta}')\leq R+\sum_{j=1}^KC_j\max_{i,j=1,\dots,K}\sum_{k=1}^{\infty}(t_k-t_{k-1})\sup_{t\in[t_{k-1},t_k],\pmb{\eta}'\in \Theta_{\pmb{\eta}}}g_{ij}(t;\pmb{\eta}'):=\Lambda_0,\\
&\sup_{\pmb{\theta}'\in \Theta}|\log\lambda_i(0;\pmb{\theta}')|\leq \max\{|\log\underline{\mu}|, ~\sup\nolimits_{\pmb{\theta}'\in \Theta}\lambda_i(0;\pmb{\theta}')\}\leq \Lambda_0+|\log\underline{\mu}|:=\Lambda_1,\\
\end{split}
\]
where $\Theta\subseteq B(0,R)$.
In addition, 
\[
\begin{split}
&\sup_{\pmb{\theta}'\in \Theta}|\lambda_i(t;\pmb{\theta}')-\hat{\lambda}_i(t;\pmb{\theta}')|\leq \sum_{j=1}^KC_j\sum_{k=1}^{\infty}(t_k-t_{k-1})\sup_{t'\in[t+t_{k-1},t+t_k],\pmb{\eta}'\in \Theta_{\pmb{\eta}}}g_{ij}(t';\pmb{\eta}').\\
\end{split}
\]

Hence it suffices to prove that $\Lambda_0$ has a finite $(3+\alpha)$-th moment for any $\alpha\in[0,1)$, and that 
\[
 \Lambda_0^{(t)}:=\sum_{j=1}^KC_j\max_{i,j=1,\dots,K}\sum_{k=1}^{\infty}(t_k-t_{k-1})\sup_{t'\in[t+t_{k-1},t+t_k],\pmb{\eta}'\in \Theta_{\pmb{\eta}}}g_{ij}(t';\pmb{\eta}')\rightarrow 0
\]
in the mean-square sense as $t\rightarrow\infty$, and is uniformly bounded for all $t\geq 0$. Here $\Lambda_0^{(0)}+R=\Lambda_0$. 

Now by Assumption \ref{assumption3} and by setting $T=0$,  
\[
\begin{split}
\tilde{C}&:=\max_{i,j=1,\dots,K}\sum_{k=1}^{\infty}(t_k-t_{k-1})\sup_{t'\in[t_{k-1},t_k],\pmb{\eta}'\in \Theta_{\pmb{\eta}}}g_{ij}(t';\pmb{\eta}')<\infty,\\
\tilde{C}^{(t)}&:=\max_{i,j=1,\dots,K}\sum_{k=1}^{\infty}(t_k-t_{k-1})\sup_{t'\in[t+t_{k-1},t+t_k],\pmb{\eta}'\in \Theta_{\pmb{\eta}}}g_{ij}(t';\pmb{\eta}')\rightarrow 0~\mbox{as $t\rightarrow\infty$},\\
\tilde{C}^{(t)}&\leq E~\mbox{for any $t\geq 0$ and some $E>0$}.
\end{split}
\]

By taking $t=0$ in Lemma \ref{uniform_increment}, we have $\mathbb{E}[|C_j|^{3+\alpha}]<\infty$ for $j=1,\dots,K$ for any $\alpha\in[0,1)$. Since $\Lambda_0\leq R+\tilde{C}\sum_{j=1}^KC_j$ and $\Lambda_0^{(t)}\leq \tilde{C}^{(t)}\sum_{j=1}^KC_j$, we see that 
\begin{itemize}
\item $\mathbb{E}[|\Lambda_0|^{3+\alpha}]<\infty$;
\item $\mathbb{E}[|\Lambda_0^{(t)}|^{3+\alpha}]\leq F$ for all $t\geq 0$ and some constant $F>0$, $\mathbb{E}[|\Lambda_0^{(t)}|^{3+\alpha}]\rightarrow 0$ as $t\rightarrow\infty$, which then imply that $\mathbb{E}[|\Lambda_0^{(t)}|^{2}]\leq F$ for all $t\geq 0$ and some constant $F>0$, $\mathbb{E}[|\Lambda_0^{(t)}|^2]\rightarrow 0$ as $t\rightarrow\infty$. 
\end{itemize}

\qed

\subsection{Proof of Lemma \ref{cont_lambda}.} To prove the almost sure continuity of $\lambda_i(t;\btheta)$ as a function of $\btheta$ for each fixed $t\geq 0$, it suffices to prove the continuity for each dimension $i$ and time $t$. Below we focus on $\lambda_i(t;\btheta)$ for a fixed $i=1,\dots,K$ and $t\geq 0$. 

By the definition of $C_i^{(t)}$ in Lemma \ref{uniform_increment}, for an arbitrary sequence $\btheta_n\rightarrow\btheta_0$, where $\btheta_n=(\bmu^n,\pmb{\eta}^n)$ and $\btheta_0=(\bmu^0,\pmb{\eta}^0)$,
\[
\begin{split}
\left|\lambda_i(t;\btheta_n)-\lambda_i(t;\btheta_0)\right|&\leq \left|\mu_i^n-\mu_i^0\right|+\sum_{j=1}^K\int_{-\infty}^t\left|g_{ij}(t-s;\pmb{\eta}^n)-g_{ij}(t-s;\pmb{\eta}^0)\right|N_j(ds)\\
&\leq \left|\mu_i^n-\mu_i^0\right|+\sum_{j=1}^KC_j^{(t)}\sum_{k=1}^{\infty}(t_k-t_{k-1})\sup_{t'\in[t_{k-1},t_k]}\left|g_{ij}(t';\pmb{\eta}^n)-g_{ij}(t';\pmb{\eta}^0)\right|. \\
\end{split}
\]

Let $L_{t_{k-1},t_k}:=\sum_{d=1}^D\sup_{t'\in[t_{k-1},t_k],\pmb{\eta}'\in\Theta_{\bigeta}}\left|\partial_{\eta_d}g_{ij}(t';\bigeta')\right|$. Then by Assumption \ref{assumption3}, 
\begin{equation}\label{lip_sum}
\sum_{k=1}^{\infty}(t_k-t_{k-1})L_{t_{k-1},t_k}<\infty.
\end{equation}

By Assumption \ref{assumption1}, $\exists\epsilon>0$ such that $B(\btheta_0,\epsilon)\subseteq\tilde{\Theta}$. Since for a sufficiently large $n$, $\btheta^n\in B(\btheta_0,\epsilon)$, by the mean-value theorem, there exists some $c\in(0,1)$, such that 
\[
\begin{split}
\sup_{t'\in[t_{k-1},t_k]}\left|g_{ij}(t';\pmb{\eta}^n)-g_{ij}(t';\pmb{\eta}^0)\right|&=\sup_{t'\in[t_{k-1},t_k]}\left|\nabla_{\bigeta} g_{ij}((1-c)\bigeta^0+c~\bigeta^n)^T(\bigeta^n-\bigeta^0)\right| \\
&\leq \sup_{t'\in[t_{k-1},t_k]}\sup_{\bigeta'\in\Theta_{\bigeta}}\|\nabla_{\bigeta}g_{ij}(\bigeta')\|_1\|\bigeta^n-\bigeta^0\|_{\infty}\leq L_{t_{k-1},t_k}\|\bigeta^n-\bigeta^0\|_{\infty},
\end{split}
\]
hence
\begin{equation}\label{bound_cont}
\left|\lambda_i(t;\btheta_n)-\lambda_i(t;\btheta_0)\right|\leq 
\left|\mu_i^n-\mu_i^0\right|+\sum_{j=1}^KC_j^{(t)}\sum_{k=1}^{\infty}(t_k-t_{k-1})L_{t_{k-1},t_k}\left\|\pmb{\eta}^n-\pmb{\eta}^0\right\|_{\infty}.
\end{equation}

By Lemma \ref{uniform_increment}, we have in particular that $C_j^{(t)}$ is a.s. finite for any $j=1,\dots,K$. Hence as $\btheta_n\rightarrow\btheta_0$, $\mu_i^n\rightarrow\mu_i^0$ and $\bigeta_n\rightarrow\bigeta_0$ as $n\rightarrow\infty$, and hence  $\lambda_i(t;\btheta_n)\rightarrow\lambda_i(t;\btheta_0)$ almost surely as $n\rightarrow\infty$ from Eqns. (\ref{lip_sum}) and (\ref{bound_cont}). 
\qed

\subsection{Proof of Lemma \ref{ergodic_xi}.} 
{Since $\pmb{N}$ is stationary, the time shift operator $S_1$ through the unit distance is measure-preserving [\cite{Vere-Jones}, Chapter 12.2]}. Moreover, {by the adaptedness}, for any $t\geq 0$, $\xi(t)$ is a measurable functional of the point process $\{N_i(s,t'],s<t'<t,i=1,\dots,K\}$. Along with the {\color{black}assumption that the underlying (true) MHP is ergodic}, the $\sigma$-algebra of invariant events under $S_1$ is trivial. Since $\mathbb{E}\left[\left|\int_0^1\xi(t)dt\right|\right]\leq \mathbb{E}[|\xi(0)|]<\infty$, applying Birkhoff's ergodic theorem (\cite{Durrett}) to $X:=\int_0^1\xi(t)dt$ and $S_1$ yields
\[
\lim_{N\rightarrow\infty}\dfrac{1}{N}\int_0^N\xi(t)dt=\lim_{N\rightarrow\infty}\dfrac{1}{N}\sum_{k=0}^{N-1}\int_{k}^{k+1}\xi(t)dt=\mathbb{E}[\xi(0)],
\]
where $N$ takes positive integer values. And since by assumption $\mathbb{E}[|\xi(0)|^2]<\infty$, we have for general $T>0$,
\[
\left|\dfrac{\int_0^T\xi(t)dt-\int_0^{\lfloor T\rfloor}\xi(t)dt}{\lfloor T\rfloor}\right|=\left|\dfrac{\int_{\lfloor T\rfloor}^T\xi(t)dt}{\lfloor T\rfloor}\right|\leq \dfrac{1}{\lfloor T\rfloor}\int_{T-1}^T|\xi(t)|dt\rightarrow 0~\text{in probability},
\]
where the last limit follows from $\mathbb{P}\left(Y(T)\geq \epsilon\right)\leq \dfrac{\mathbb{E}[Y(T)]}{\epsilon}=\dfrac{\mathbb{E}|\xi(0)|}{\epsilon \lfloor T\rfloor}\rightarrow 0$ as $T\rightarrow\infty$, according to the definition $Y(T):=\frac{1}{\lfloor T\rfloor}\int_{T-1}^T|\xi(t)|dt$.


Together with the fact that 
$\lim_{T\rightarrow\infty}\dfrac{1}{T}\int_0^T\xi(t)dt=\lim_{T\rightarrow\infty}\dfrac{1}{\lfloor T\rfloor}\int_0^{T}\xi(t)dt$ a.s., we conclude that
\[
\lim_{T\rightarrow\infty}\dfrac{1}{T}\int_0^T\xi(t)dt=\lim_{T\rightarrow\infty}\dfrac{1}{\lfloor T\rfloor}\int_0^{\lfloor T\rfloor}\xi(t)dt=\mathbb{E}[\xi(0)]~\text{in probability,}
\]
where $T$ takes general positive real values. This proves Eqn. (\ref{ergodic_xi_1}).

As for Eqn. (\ref{ergodic_xi_2}), consider 
\begin{equation*}
\eta_i(T):=\int_0^T \xi(t)\frac{N_i(dt)}{\lambda_i(t;\btheta^\star)}-\int_0^T \xi(t)dt=\int_0^T\dfrac{\xi(t)}{\lambda_i(t;\btheta^\star)}M_i(dt),
\end{equation*}

Then since $\xi(t)$ is stationary and has finite second-order moments, we have
\[
\int_0^T\mathbb{E}\left[\left(\dfrac{\xi(t)}{\lambda_i(t;\btheta^\star)}\right)^2\lambda_i(t;\btheta^\star)\right]dt
\leq \dfrac{T}{\underline{\mu}}\mathbb{E}\left[\xi(0)^2\right]<\infty.
\]
Hence applying Proposition \ref{dN2dt_ogata} gives
\[
\begin{split}
\mathbb{E}\left[\int_0^T\left(\dfrac{\xi(t)}{\lambda_i(t;\btheta^\star)}\right)^2N_i(dt)\right]
&=\mathbb{E}\left[\int_0^T\left(\dfrac{\xi(t)}{\lambda_i(t;\btheta^\star)}\right)^2\lambda_i(t;\btheta^\star)dt\right]
<\infty.
\end{split}
\]

By Proposition \ref{dN2dt}, 
\[
\begin{split}
\mathbb{E}\left[\eta_i(T)^2\right]
&=\mathbb{E}\left[\int_0^T\left(\dfrac{\xi(t)}{\lambda_i(t;\btheta^\star)}\right)^2N_i(dt)\right]\leq\dfrac{T}{\underline{\mu}}\mathbb{E}\left[\xi(0)^2\right],
\end{split}
\]
from which
\[
\lim_{T\rightarrow\infty}\mathbb{E}\left[\left(\dfrac{1}{T}\eta_i(T)\right)^2\right]\leq\lim_{T\rightarrow\infty}\dfrac{1}{T\underline{\mu}}\mathbb{E}\left[\xi(0)^2\right]=0.
\]

Since convergence in expectation implies convergence in probability, $\lim_{T\rightarrow\infty}\eta_i(T)/T=0$ in probability. This, together with Eqn. (\ref{ergodic_xi_1}), implies Eqn. (\ref{ergodic_xi_2}).
\qed

\subsection{Proof of Lemma \ref{stat_moment_inflambda}.}

The left continuity of $\xi_{U,i}^{(l)}(t)$ ($l=1,2,3$) is directly implied from the left continuity of $g_{ij}$ in $t\geq 0$. 
Now we prove the bounds related to $\xi_{U,i}^{(l)}(t)$, $l=1,\dots,3$. 
First, by Lemma \ref{newlemma_consist}
\[
\mathbb{E}\left[\left|\xi_{U,i}^{(1)}(0)\right|^2\right]\leq \mathbb{E}\left[(2\Lambda_0)^2\right]=4\mathbb{E}[\Lambda_0^2]<\infty.
\]

Secondly, since $\log x\leq x-1$ for any $x\geq 0$,
\[
-\xi_{U,i}^{(2)}(0)=\lambda_i(0;\btheta^\star)\log\left(\frac{\sup_{\theta\in U}\lambda_i(0;\btheta)}{\lambda_i(0;\btheta^{\star})}\right)\leq \sup_{\theta\in U}\lambda_i(0;\btheta)-\lambda_i(0;\btheta^\star).
\] 
Meanwhile, we also have
\[
\begin{split}
-\xi_{U,i}^{(2)}(0)=\lambda_i(0;\btheta^\star)\log\left(\frac{\sup_{\btheta\in U}\lambda_i(0;\btheta)}{\lambda_i(0;\btheta^{\star})}\right)&=\lambda_i(0;\pmb{\btheta}^\star)\log(\sup\nolimits_{\theta\in U}\lambda_i(0;\btheta))-\lambda_i(0;\btheta^\star)\log\lambda_i(0;\btheta^\star)\\
&\geq \lambda_i(0;\btheta^\star)\log\underline{\mu}-\lambda_i(0;\btheta^\star)\log\lambda_i(0;\btheta^\star).
\end{split}
\] 

Since $x\log x=O(x^{1+\alpha})$ for any $\alpha>0$, we see that there exists some constant $c>0$ such that 
$x\log x\leq x^{1+\alpha}~\mbox{for all $x\geq c$}$. Hence for any $x\geq \underline{\mu}$, we have by choosing $\alpha=1/2$,
\[
x\log x\leq \max\{\underline{\mu}\left|\log\underline{\mu}\right|,~c|\log c|,~x^{3/2}\}.
\]
Now replacing $x$ with $\lambda_i(0;\btheta^\star)$, and noticing that $\lambda_i(0;\btheta^\star)\geq \underline{\mu}$, then again by Lemma \ref{newlemma_consist}
\[
\mathbb{E}\left[\left|\xi_{U,i}^{(2)}(0)\right|^2\right]\leq\mathbb{E}\left[\left(2\Lambda_0+\left|\log\underline{\mu}\right|\Lambda_0+\underline{\mu}\left|\log\underline{\mu}\right|+c|\log c|+\Lambda_0^{3/2}\right)^2\right]<\infty.
\]

Finally, by Lemma \ref{approx} and stationarity of $\lambda_i(t;\btheta^\star)$ 
\[
\begin{split}
\int_0^T\mathbb{E}&\left[\left|\xi_{U,i}^{(3)}(t)\right|\lambda_i(t;\btheta^\star)\right]dt\leq\int_0^T\sqrt{\mathbb{E}[\lambda_i^2(t;\btheta^\star)]\mathbb{E}\left[\sup_{\btheta\in U}\left|\lambda_i(t;\btheta)-\hat{\lambda}_i(t;\btheta)\right|^2\right]}dt\\
&=\sqrt{\mathbb{E}[\lambda_i^2(0;\btheta^\star)]}\int_0^T\sqrt{\mathbb{E}\left[\sup_{\btheta\in U}\left|\lambda_i(t;\btheta)-\hat{\lambda}_i(t;\btheta)\right|^2\right]}dt\leq T\sqrt{F\mathbb{E}\lambda_i^2(0;\btheta^\star)}<\infty.
\end{split}
\]
\qed

\section{Proofs for Section\ref{optimization}}\label{proofs_opt}

\subsection{Proof of Proposition \ref{glb_conv_set_crit}.}
The proof of $L\subseteq S$ is a direct combination of Propositions 4.3 and 4.5 in \cite{iPALM}. Furthermore, notice that by the definition of algorithm iPALM, all iterations are feasible, and hence bounded as $\Theta$ is compact. By Weierstrass theorem, the limit point set $L$ is non-empty. Since the limit point set is always closed, $L$ is compact. Finally, by the proof of Proposition 4.5 in \cite{iPALM}, we know that $\|\btheta^{k+1}-\btheta^k\|_2\rightarrow 0$ as $k\rightarrow\infty$.
Hence $L$ is also connected: otherwise by compactness, $L$ can be partitioned into two disjoint sets $L_1$ and $L_2$ with a positive distance, contradicting with the zero limit of $\|\btheta^{k+1}-\btheta^k\|_2$. 
\qed

\subsection{Proof of Lemma \ref{fix=crit}.}
If $(\btheta',\btheta)$ is a fixed point of $H^{iPALM}$, then by definition $\btheta=\btheta'$, and it is the limit point of a iteration sequence generated by iPALM with $\btheta^0=\btheta=\btheta'$. Hence by Proposition \ref{glb_conv_set_crit}, it is also a stationary point of $F$. 

Conversely, if $\btheta=(\bmu,\balpha,\bbeta)$ is a stationary point of $F$, then for $\btheta'=\btheta$, we have 
\[
\pmb{0}\in-\nabla_{\bmu,\balpha}\hat{L}_T^{reg}({\bmu,\balpha}) + N_A({\bmu,\balpha}),\quad \pmb{0}\in-\nabla_{\bbeta}\hat{L}_T^{reg}({\bbeta}) + N_B({\bbeta}),
\]
and since $A$ and $B$ are nonempty closed and convex, by the properties of proximal (gradient) mappings (among which projections are special cases) (\textit{e.g.}, \cite{Proximal}), we have 
 \[
\left[ \begin{array}{l}
 \bmu\\
 \balpha
 \end{array}\right]=\Pi_A\left(
 \left[ \begin{array}{l}
 \bmu\\
 \balpha
 \end{array}\right]
 +\tau_1^k\nabla_{\bmu,\balpha}\hat{L}_T^{reg}\left(\bmu,\balpha,\bbeta\right) \right),\quad
 \bbeta=\Pi_B\left(\bbeta+\tau_2^k\nabla_{\bbeta}\hat{L}_T^{reg}(\bmu,\balpha,\bbeta)\right).
 \]
 Hence $(\btheta,\btheta)$ is a fixed point of $H^{iPALM}$.
 \qed

\subsection{Proof of Theorem \ref{AA-iPALM-convergence}.}
We begin by noticing that $\pmb{u}^{k+1}$ either equals $\tilde{\pmb{u}}^{k+1}=\pmb{u}^k-H_k(\pmb{u}^k-\hat{\pmb{u}}^{k+1})$ or
$\hat{\pmb{u}}^{k+1}=H^{iPALM}(\pmb{u}^k)$, depending on whether the conditions in line 15 of Algorithm \ref{alg:AA-iPALM} are satisfied or not. We partition the iteration counts into two subsets accordingly, with $K_{AA}=\{k_0,k_1,\dots\}$ containing those iterations that pass line 15, while $K_{iPALM}=\{l_0,l_1,\dots\}$ being the rest that go to line 16. 

\paragraph{Step 1.} For $k_i\in K_{AA}$ ($i\geq 0$), by the fourth inequality in line 15 of Algorithm \ref{alg:AA-iPALM}, the fact that $\epsilon\in(0,1)$ and that $(\bmu^{k_i+1},\balpha^{k_i+1},\bbeta^{k_i+1})=(\tilde{\bmu}^{k_i+1},\tilde{\balpha}^{k_i+1},\tilde{\bbeta}^{k_i+1})$, we have
\begin{equation*}
\begin{split}
\rho_1&\left(\|{\bmu}^{k_i+1}-\bmu^{k_i}\|_2^2+\|{\balpha}^{k_i+1}-\balpha^{k_i}\|_2^2+\|{\bbeta}^{k_i+1}-\bbeta^{k_i}\|_2^2+\|\bmu^{k_i}-\bmu^{k_i-1}\|_2^2+\|\balpha^{k_i}-\balpha^{k_i-1}\|_2^2+\|\bbeta^{k_i}-\bbeta^{k_i-1}\|_2^2\right)\\
\leq &\frac{\epsilon\delta}{2}\left(\|\tilde{\bmu}^{k_i+1}-\bmu^{k_i}\|_2^2+\|\tilde{\balpha}^{k_i+1}-\balpha^{k_i}\|_2^2+\|\tilde{\bbeta}^{k_i+1}-\bbeta^{k_i}\|_2^2\right)\\
&+\frac{\delta_1}{2}\left(\|\bmu^{k_i}-\bmu^{k_i-1}\|_2^2+\|\balpha^{k_i}-\balpha^{k_i-1}\|_2^2\right)+\frac{\delta_2}{2}\|\bbeta^{k_i}-\bbeta^{k_i-1}\|_2^2\\
\leq &\hat{L}_T^{reg}(\tilde{\bmu}^{k_i+1},\tilde{\balpha}^{k_i+1},\tilde{\bbeta}^{k_i+1})-\hat{L}_T^{reg}(\bmu^{k_i},\balpha^{k_i},\bbeta^{k_i})- \dfrac{\delta}{2}(\|\tilde{\bmu}^{k_i+1}-\bmu^{k_i}\|_2^2+\|\tilde{\balpha}^{k_i+1}-\balpha^{k_i}\|_2^2+\|\tilde{\bbeta}^{k_i+1}-\bbeta^{k_i}\|_2^2)\\
&+\frac{\delta_1}{2}\left(\|\bmu^{k_i}-\bmu^{k_i-1}\|_2^2+\|\balpha^{k_i}-\balpha^{k_i-1}\|_2^2\right)+\frac{\delta_2}{2}\|\bbeta^{k_i}-\bbeta^{k_i-1}\|_2^2\\
\leq & \hat{L}_T^{reg}({\bmu}^{k_i+1},{\balpha}^{k_i+1},{\bbeta}^{k_i+1})-\hat{L}_T^{reg}(\bmu^{k_i},\balpha^{k_i},\bbeta^{k_i}) \\
&+ \frac{\delta_1}{2}\left(\|(\bmu^{k_i}-\bmu^{k_i-1},\balpha^{k_i}-\balpha^{k_i-1})\|_2^2-\|({\bmu}^{k_i+1}-\bmu^{k_i},{\balpha}^{k_i+1}-\balpha^{k_i})\|_2^2\right)+\frac{\delta_2}{2}\left(\|\bbeta^{k_i}-\bbeta^{k_i-1}\|_2^2-\|{\bbeta}^{k_i+1}-\bbeta^{k_i}\|_2^2\right).
\end{split}
\end{equation*}
Meanwhile, Proposition \ref{ipalm_bounds} implies that the same inequality also holds for $l_i\in K_{iPALM}$. Hence for any $k\geq 0$,
\begin{equation*}
\begin{split}
\rho_1&\left(\|{\bmu}^{k+1}-\bmu^k\|_2^2+\|{\balpha}^{k+1}-\balpha^k\|_2^2+\|{\bbeta}^{k+1}-\bbeta^k\|_2^2+\|\bmu^k-\bmu^{k-1}\|_2^2+\|\balpha^k-\balpha^{k-1}\|_2^2+\|\bbeta^k-\bbeta^{k-1}\|_2^2\right)\\
\leq & \hat{L}_T^{reg}({\bmu}^{k+1},{\balpha}^{k+1},{\bbeta}^{k+1})-\hat{L}_T^{reg}(\bmu^k,\balpha^k,\bbeta^k) + \frac{\delta_1}{2}\left(\|(\bmu^k-\bmu^{k-1},\balpha^k-\balpha^{k-1})\|_2^2-\|({\bmu}^{k+1}-\bmu^k,{\balpha}^{k+1}-\balpha^k)\|_2^2\right)\\
&+\frac{\delta_2}{2}\left(\|\bbeta^k-\bbeta^{k-1}\|_2^2-\|{\bbeta}^{k+1}-\bbeta^k\|_2^2\right).
\end{split}
\end{equation*}
 
\paragraph{Step 2.} For $k_i\in K_{AA}$ ($i\geq 0$), by inequality (\ref{ipalm_bound2}), 
\begin{equation*}
\begin{split}
&\left\|\left(\nabla_{\bmu,\balpha}\hat{L}_T^{reg}({\bmu}^{k_i},{\balpha}^{k_i},{\bbeta}^{k_i}), ~\nabla_{\bbeta}\hat{L}_T^{reg}(\hat{\bmu}^{k_i+1},\hat{\balpha}^{k_i+1},{\bbeta}^{k_i})\right)-\nabla\hat{L}_T^{reg}(\hat{\bmu}^{k_i+1},\hat{\balpha}^{k_i+1},\hat{\bbeta}^{k_i+1})\right.\\
&+\left(\frac{1}{\tau_1}\left((\bmu^{k_i}-\hat{\bmu}^{k_i+1},\balpha^{k_i}-\hat{\balpha}^{k_i+1})+\gamma_1(\bmu^{k_i}-\bmu^{k_i-1},\balpha^{k_i}-\balpha^{k_i-1})\right),\frac{1}{\tau_2}\left(\bbeta^{k_i}-\hat{\bbeta}^{k_i+1}+\gamma_2(\bbeta^{k_i}-\bbeta^{k_i-1})\right)\right)\\
&\left. + \left(\delta_1(\hat{\bmu}^{k_i+1}-\bmu^{k_i},\hat{\balpha}^{k_i+1}-\balpha^{k_i}),\delta_2(\hat{\bbeta}^{k_i+1}-\bbeta^{k_i})\right)\right\|_2\\
&\leq \rho_2\|\hat{\pmb{u}}^{k_i+1}-\pmb{u}^{k_i}\|_2.
\end{split}
\end{equation*}

Now by Proposition \ref{Hkbounds}, the definition of $M$ and the fact that $(\bmu^{k_i+1},\balpha^{k_i+1},\bbeta^{k_i+1})=(\tilde{\bmu}^{k_i+1},\tilde{\balpha}^{k_i+1},\tilde{\bbeta}^{k_i+1})$ (due to the second condition in line 15 of Algorithm \ref{alg:AA-iPALM}), 
\begin{equation*}
\begin{split}
&\left\|\nabla\hat{L}_T^{reg}({\bmu}^{k_i+1}, {\balpha}^{k_i+1}, {\bbeta}^{k_i+1})-
\nabla\hat{L}_T^{reg}(\hat{\bmu}^{k_i+1},\hat{\balpha}^{k_i+1},\hat{\bbeta}^{k_i+1})\right\|_2\\
&\leq M\left\|\left(\tilde{\bmu}^{k_i+1}-\hat{\bmu}^{k_i+1},\tilde{\balpha}^{k_i+1}-\hat{\balpha}^{k_i+1},\tilde{\bbeta}^{k_i+1}-\hat{\bbeta}^{k_i+1}\right)\right\|_2\\
&\leq M\|\tilde{\pmb{u}}^{k_i+1}-\hat{\pmb{u}}^{k_i+1}\|_2=M\|(I-H_{k_i})(\pmb{u}^{k_i}-\hat{\pmb{u}}^{k_i+1})\|_2
\leq M(1+\sigma_H^+)\|\pmb{u}^{k_i}-\hat{\pmb{u}}^{k_i+1}\|_2.
\end{split}
\end{equation*}

Similarly,
\begin{equation*}
\begin{split}
&\|({\bmu}^{k_i+1}-\bmu^{k_i},{\balpha}^{k_i+1}-\balpha^{k_i}, {\bbeta}^{k_i+1}-\bbeta^{k_i})-(\hat{\bmu}^{k_i+1}-\bmu^{k_i},\hat{\balpha}^{k_i+1}-\balpha^{k_i}, \hat{\bbeta}^{k_i+1}-\bbeta^{k_i})\|_2\\
&\leq (1+\sigma_H^+)\|\pmb{u}^{k_i}-\hat{\pmb{u}}^{k_i+1}\|_2.
\end{split}
\end{equation*}

In addition, by the first condition in line 15 of Algorithm \ref{alg:AA-iPALM},
\begin{equation*}
\begin{split}
\left\|\left(\nabla_{\bmu,\balpha}\hat{L}_T^{reg}({\bmu}^{k_i},{\balpha}^{k_i},{\bbeta}^{k_i}),~\nabla_{\bbeta}\hat{L}_T^{reg}({\bmu}^{k_i},{\balpha}^{k_i},{\bbeta}^{k_i})\right)\right\|_2\leq C_1\left\|\hat{\pmb{u}}^{k_i+1}-\pmb{u}^{k_i}\right\|_2,
\end{split}
\end{equation*}
and hence by $\|\nabla_{\bbeta}\hat{L}_T^{reg}(\hat{\bmu}^{k_i+1},\hat{\balpha}^{k_i+1},{\bbeta}^{k_i})-\nabla_{\bbeta}\hat{L}_T^{reg}({\bmu}^{k_i},{\balpha}^{k_i},{\bbeta}^{k_i})\|_2\leq L_2\left\|\hat{\pmb{u}}^{k_i+1}-\pmb{u}^{k_i}\right\|_2$, 
\begin{equation*}
\begin{split}
\left\|\left(\nabla_{\bmu,\balpha}\hat{L}_T^{reg}({\bmu}^{k_i},{\balpha}^{k_i},{\bbeta}^{k_i}),~\nabla_{\bbeta}\hat{L}_T^{reg}(\hat{\bmu}^{k_i+1},\hat{\balpha}^{k_i+1},{\bbeta}^{k_i})\right)\right\|_2\leq \sqrt{2C_1^2+2L_2^2}\left\|\hat{\pmb{u}}^{k_i+1}-\pmb{u}^{k_i}\right\|_2.
\end{split}
\end{equation*}

Since $\pmb{0}\in N_A(\tilde{\bmu}^{k+1},\tilde{\balpha}^{k+1})$, $\pmb{0}\in N_B(\tilde{\bbeta}^{k+1})$, we see that $\exists$ $\pmb{g}^{k_i+1}\in \partial F({\bmu}^{k_i+1},{\balpha}^{k_i+1},{\bbeta}^{k_i+1})+\left(\delta_1({\bmu}^{k_i+1}-\bmu^{k_i},\right.$ $\left.{\balpha}^{k_i+1}-\balpha^{k_i}), \delta_2(\bbeta^{k_i+1}-\bbeta^{k_i})\right)$, such that 
\begin{equation*}
\begin{split}
\|\pmb{g}^{k_i+1}\|_2\leq \left(\rho_2+M(1+\sigma_H^+)+\sqrt{2C_1^2+2L_2^2}+\left(\delta+\frac{1}{\tau}\right)(1+\sigma_H^+)+\frac{\gamma}{\tau}\right)\|\hat{\pmb{u}}^{k_i+1}-\pmb{u}^{k_i}\|_2,
\end{split}
\end{equation*}
where $\tau:=\min\{\tau_1,\tau_2\}$, $\gamma:=\max\{\gamma_1,\gamma_2\}$.

Denoting $\tilde{\pmb{u}}^{k_i+1}_2=(\tilde{\bmu}^{k_i},\tilde{\balpha}^{k_i},\tilde{\bbeta}^{k_i})$, then by the third inequality in line 15 of Algorithm \ref{alg:AA-iPALM} and again Proposition \ref{Hkbounds},
\begin{equation*}
\begin{split}
\|\pmb{u}^{k_i}-\hat{\pmb{u}}^{k_i+1}\|_2^2&\leq \|H_{k_i}^{-1}\|_2^2 \|H_{k_i}(\pmb{u}^{k_i}-\hat{\pmb{u}}^{k_i+1})\|_2^2=
\|H_{k_i}^{-1}\|_2^2\|\tilde{\pmb{u}}^{k_i+1}-\pmb{u}^{k_i}\|_2^2\\
\leq& (\sigma_H^-)^2\left(\|{\bmu}^{k_i+1}-\bmu^{k_i}\|_2^2+\|{\balpha}^{k_i+1}-\balpha^{k_i}\|_2^2+\|{\bbeta}^{k_i+1}-\bbeta^{k_i}\|_2^2 + \|\tilde{\pmb{u}}^{k_i+1}_2-{\pmb{u}}^{k_i}_2\|_2^2\right)\\
&\leq (C_2\sigma_H^-)^2\left(\|{\bmu}^{k_i+1}-\bmu^{k_i}\|_2^2+\|{\balpha}^{k_i+1}-\balpha^{k_i}\|_2^2+\|{\bbeta}^{k_i+1}-\bbeta^{k_i}\|_2^2\right.\\
&\left.+\|{\bmu}^{k_i}-\bmu^{k_i-1}\|_2^2+\|{\balpha}^{k_i}-\balpha^{k_i-1}\|_2^2+\|{\bbeta}^{k_i}-\bbeta^{k_i-1}\|_2^2\right).
\end{split}
\end{equation*}

Hence
\begin{equation*}
\begin{split}
\|\pmb{g}^{k_i+1}\|_2^2&\leq 
\left(C_2\sigma_H^-\left(\rho_2+\sqrt{2C_1^2+2L_2^2}+\left(M+\delta+\frac{1}{\tau}\right)(1+\sigma_H^+)+\frac{\gamma}{\tau}\right)\right)^2
\left(\|{\bmu}^{k_i+1}-\bmu^{k_i}\|_2^2\right.\\
&\left.+\|{\balpha}^{k_i+1}-\balpha^{k_i}\|_2^2+\|{\bbeta}^{k_i+1}-\bbeta^{k_i}\|_2^2+\|{\bmu}^{k_i}-\bmu^{k_i-1}\|_2^2+\|{\balpha}^{k_i}-\balpha^{k_i-1}\|_2^2+\|{\bbeta}^{k_i}-\bbeta^{k_i-1}\|_2^2\right).
\end{split}
\end{equation*}

Moreover, for $l_i\in K_{iPALM}$ ($i\geq 0$), by inequality (\ref{ipalm_bound2}), we immediately see that the same inequality holds with a smaller constant $\rho_2^2$ (since $C_2,~\sigma_H^-\geq 1$ by definition). 
Hence for any $k\geq 0$, $\exists$ $\pmb{g}^{k+1}\in \partial F({\bmu}^{k+1},{\balpha}^{k+1},{\bbeta}^{k+1})+\left(\delta_1({\bmu}^{k+1}-\bmu^{k},{\balpha}^{k+1}-\balpha^{k}), \delta_2(\bbeta^{k+1}-\bbeta^{k})\right)$, such that 
\begin{equation*}
\begin{split}
\|\pmb{g}^{k+1}\|_2^2&\leq 
\left(C_2\sigma_H^-\left(\rho_2+\sqrt{2C_1^2+2L_2^2}+\left(M+\delta+\frac{1}{\tau}\right)(1+\sigma_H^+)+\frac{\gamma}{\tau}\right)\right)^2
\left(\|{\bmu}^{k+1}-\bmu^{k}\|_2^2\right.\\
&\left.+\|{\balpha}^{k+1}-\balpha^{k}\|_2^2+\|{\bbeta}^{k+1}-\bbeta^{k}\|_2^2+\|{\bmu}^{k}-\bmu^{k-1}\|_2^2+\|{\balpha}^{k}-\balpha^{k-1}\|_2^2+\|{\bbeta}^{k}-\bbeta^{k-1}\|_2^2\right).
\end{split}
\end{equation*} 

\paragraph{Step 3.}  Finally, combining steps 1 and 2,  and recalling that $\Theta\in B(0,R)$ (notice that all iterations are in $\Theta$ by the safeguarding conditions), we see that for any $K\geq 1$,
\begin{equation*}
\begin{split}
\sum_{k=0}^{K-1}&\|\pmb{g}^{k+1}\|_2^2\leq \left(C_2\sigma_H^-\left(\rho_2+\sqrt{2C_1^2+2L_2^2}+\left(M+\delta+\frac{1}{\tau}\right)(1+\sigma_H^+)+\frac{\gamma}{\tau}\right)\right)^2 \sum_{k=0}^{K-1}\left(\|{\bmu}^{k+1}-\bmu^{k}\|_2^2\right.\\
&\left.+\|{\balpha}^{k+1}-\balpha^{k}\|_2^2+\|{\bbeta}^{k+1}-\bbeta^{k}\|_2^2+\|{\bmu}^{k}-\bmu^{k-1}\|_2^2+\|{\balpha}^{k}-\balpha^{k-1}\|_2^2+\|{\bbeta}^{k}-\bbeta^{k-1}\|_2^2\right)\\
\leq & \dfrac{\left(C_2\sigma_H^-\left(\rho_2+\sqrt{2C_1^2+2L_2^2}+\left(M+\delta+\frac{1}{\tau}\right)(1+\sigma_H^+)+\frac{\gamma}{\tau}\right)\right)^2}{\rho_1}\left(\hat{L}_T^{reg}(\bmu^K,\balpha^K,\bbeta^K)-\hat{L}_T^{reg}(\bmu^0,\balpha^0,\bbeta^0)\right.\\
&\left.+\frac{\delta_1}{2}\|(\bmu^{K-1}-\bmu^{K-2},\balpha^{K-1}-\balpha^{K-2})\|_2^2+\frac{\delta_2}{2}\|\bbeta^{K-1}-\bbeta^{K-2}\|_2^2\right)\\
\leq & \dfrac{\left(C_2\sigma_H^-\left(\rho_2+\sqrt{2C_1^2+2L_2^2}+\left(M+\delta+\frac{1}{\tau}\right)(1+\sigma_H^+)+\frac{\gamma}{\tau}\right)\right)^2}{\rho_1}\\
&\times \left(\sup_{(\bmu,\balpha,\bbeta)\in \Theta}\hat{L}_T^{reg}(\bmu,\balpha,\bbeta)-\hat{L}_T^{reg}(\bmu^0,\balpha^0,\bbeta^0)+2{\delta}R^2\right).
\end{split}
\end{equation*}

This implies that $\|{\bmu}^{k+1}-\bmu^{k}\|_2^2+\|{\balpha}^{k+1}-\balpha^{k}\|_2^2+\|{\bbeta}^{k+1}-\bbeta^{k}\|_2^2+\|{\bmu}^{k}-\bmu^{k-1}\|_2^2+\|{\balpha}^{k}-\balpha^{k-1}\|_2^2+\|{\bbeta}^{k}-\bbeta^{k-1}\|_2^2\rightarrow 0$ as $k\rightarrow\infty$, which then implies that $\|\pmb{g}^{k+1}\|_2\rightarrow0$ as well as $(\bmu^{k+1}-\bmu^k,\balpha^{k+1}-\balpha^k,\bbeta^{k+1}-\bbeta^k)\rightarrow 0$ as $k\rightarrow\infty$ (by step 2). Hence for any limit point $(\bmu^\star,\balpha^\star,\bbeta^\star)$ of $(\bmu^k,\balpha^k,\bbeta^k)$, by the closedness of the sub-differential mapping, $\pmb{0}\in (\partial_{\bmu,\balpha}F(\bmu^\star,\balpha^\star,\bbeta^\star),\partial_{\bbeta}F(\bmu^\star,\balpha^\star,\bbeta^\star))$, \textit{i.e.}, $(\bmu^\star,\balpha^\star,\bbeta^\star)$ is a stationary point of $F$. Notice here $\partial_{\bmu,\balpha}F(\bmu^\star,\balpha^\star,\bbeta^\star)=N_A(\bmu^\star,\balpha^\star)-\nabla_{\bmu,\balpha}\hat{L}_T^{reg}(\bmu^\star,\balpha^\star,\bbeta^\star)$ and $\partial_{\bbeta}F(\bmu^\star,\balpha^\star,\bbeta^\star)=N_B(\bbeta^\star)-\nabla_{\bbeta}\hat{L}_T^{reg}(\bmu^\star,\balpha^\star,\bbeta^\star)$. The fact that the limit point set is nonempty, compact and connected follows exactly the same proof of Proposition \ref{glb_conv_set_crit} (\textit{e.g.}, Appendix \ref{proofs_opt}). \qed

\subsection{Proof of Theorem \ref{AA-iPALM-complexity}.}
Recall in the proof of Theorem \ref{AA-iPALM-convergence} that
\begin{equation*}
\begin{split}
\sum_{k=0}^{K-1}\|\pmb{g}^{k+1}\|_2^2\leq&  \dfrac{\left(C_2\sigma_H^-\left(\rho_2+\sqrt{2C_1^2+2L_2^2}+\left(M+\delta+\frac{1}{\tau}\right)(1+\sigma_H^+)+\frac{\gamma}{\tau}\right)\right)^2}{\rho_1}\\
&\times\left(\sup_{(\bmu,\balpha,\bbeta)\in \Theta}\hat{L}_T^{reg}(\bmu,\balpha,\bbeta)-\hat{L}_T^{reg}(\bmu^0,\balpha^0,\bbeta^0)+2{\delta}R^2\right).
\end{split}
\end{equation*}
Recall also
\begin{equation*}
\begin{split}
\sum_{k=0}^{K-1}&\left(\left\|\delta_1(\bmu^{k+1}-\bmu^k,\balpha^{k+1}-\balpha^k)\right\|_2^2+\left\|\delta_2(\bbeta^{k+1}-\bbeta^k)\right\|_2^2\right)\\
\leq &\delta^2\sum_{k=0}^{K-1}\left(\|{\bmu}^{k+1}-\bmu^{k}\|_2^2+\|{\balpha}^{k+1}-\balpha^{k}\|_2^2+\|{\bbeta}^{k+1}-\bbeta^{k}\|_2^2\right.\\
&\left.+\|{\bmu}^{k}-\bmu^{k-1}\|_2^2+\|{\balpha}^{k}-\balpha^{k-1}\|_2^2+\|{\bbeta}^{k}-\bbeta^{k-1}\|_2^2\right)\\
\leq&\dfrac{\delta^2}{\rho_1}\left(\sup_{(\bmu,\balpha,\bbeta)\in \Theta}\hat{L}_T^{reg}(\bmu,\balpha,\bbeta)-\hat{L}_T^{reg}(\bmu^0,\balpha^0,\bbeta^0)+2{\delta}R^2\right),
\end{split}
\end{equation*}
we conclude from $\|a+b\|_2^2\leq 2\|a\|_2^2+2\|b\|_2^2$ that
\begin{equation*}
\begin{split}
&\sum_{k=0}^{K-1}~\textbf{dist}(\pmb{0},\partial_{\bmu,\balpha,\bbeta}F(\bmu^k,\balpha^k,\bbeta^k))^2\\
&\leq \dfrac{2\delta^2+2\left(C_2\sigma_H^-\left(\rho_2+\sqrt{2C_1^2+2L_2^2}+\left(M+\delta+\frac{1}{\tau}\right)(1+\sigma_H^+)+\frac{\gamma}{\tau}\right)\right)^2}{\rho_1}\\
&\times\left(\sup_{(\bmu,\balpha,\bbeta)\in \Theta}\hat{L}_T^{reg}(\bmu,\balpha,\bbeta)-\hat{L}_T^{reg}(\bmu^0,\balpha^0,\bbeta^0)+2{\delta}R^2\right).
\end{split}
\end{equation*}
\qed

\end{document}